\documentclass[american,preprint,3p]{elsarticle}
\pdfoutput=1
\biboptions{longnamesfirst, sort&compress}
\usepackage{babel}
\usepackage[latin1]{inputenc}
\usepackage[T1]{fontenc}
\usepackage{lmodern}
\usepackage[babel=true]{microtype}
\usepackage{amsmath,amssymb,amsthm,booktabs,color,subcaption}
\usepackage{graphicx}
\usepackage[export]{adjustbox}
\usepackage[pdftex, bookmarks]{hyperref}
\usepackage[all]{hypcap}

\usepackage{cleveref}
\crefname{equation}{}{}
\crefname{table}{Table}{Tables}
\crefname{figure}{Fig.}{Figs.}
\Crefname{figure}{Figure}{Figures}
\crefname{section}{section}{sections}
\Crefname{section}{Section}{Sections}

\usepackage{siunitx}
\DeclareMathOperator{\diam}{diam}

\DeclareMathOperator{\tr}{tr}

\newcommand{\R}{\mathbb{R}}
\newcommand{\T}{\mathcal{T}}

\providecommand{\FigPath}{.}
\date{2009}

\begin{document}

\begin{frontmatter}

\journal{Journal of Computational Physics}
\title{A new anisotropic mesh adaptation method based upon hierarchical
a~posteriori error estimates}

\author[addressKU]{Weizhang Huang}
\ead{huang@math.ku.edu}
\author[addressTUD]{Lennard Kamenski}
\ead{kamenski@mathematik.tu-darmstadt.de}
\author[addressTUD,addressCSI]{Jens Lang}
\ead{lang@mathematik.tu-darmstadt.de}

\address[addressKU]{Department of Mathematics, The University of Kansas, 
   405 Snow Hall, Lawrence, KS~66045, USA}
\address[addressTUD]{Department of Mathematics, 
   Technische Universit{\"a}t Darmstadt,
   Dolivostr.~15, D-64293 Darmstadt, Germany}
\address[addressCSI]{Center of Smart Interfaces, 
   Technische Universit{\"a}t Darmstadt,
   Petersenstr.~32, D-64287, Germany}

\begin{abstract}
   A new anisotropic mesh adaptation strategy for finite element solution of elliptic differential equations is presented. It generates anisotropic adaptive meshes as quasi-uniform ones in some metric space, with the metric tensor being computed based on hierarchical a posteriori error estimates. A global hierarchical error estimate is employed in this study to obtain reliable directional information of the solution. Instead of solving the global error problem exactly, which is costly in general, we solve it iteratively using the symmetric Gauß--Seidel (GS) method. Numerical results show that a few GS iterations are sufficient for obtaining a reasonably good approximation to the error for use in anisotropic mesh adaptation. The new method is compared with several strategies using local error estimators or recovered Hessians. Numerical results are presented for a selection of test examples and a mathematical model for heat conduction in a thermal battery with large orthotropic jumps in the material coefficients.
\end{abstract}

\begin{keyword} 
mesh adaptation \sep anisotropic mesh \sep finite elements
   \sep a~posteriori estimators
\MSC 65N50 \sep 65N30 \sep 65N15
\\[0.6\baselineskip]%
{\small{This is a preprint of a contibution published by Elsevier Inc.~in
\emph{J.~Comput.~Phys.}, 229(6) (2010), pp.~2179--2198.}}%
\\[0.4\baselineskip]%
{\footnotesize{}%
   \begin{tabular}{@{}p{0.14\linewidth}@{}p{0.85\linewidth}@{}}
      \includegraphics[height=2.2\baselineskip, valign=t]{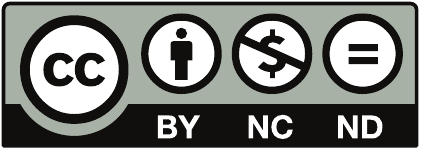}
   &
   \copyright~2009. Licensed under CC-BY-NC-ND~4.0 (\url{https://creativecommons.org/licenses/by-nc-nd/4.0}).
   \vspace{0.3\baselineskip}
   \newline{}
   The final version is available online at \url{https://dx.doi.org/10.1016/j.jcp.2009.11.029}.
   \end{tabular}%
}%
\end{keyword}
\end{frontmatter}

\section{Introduction}
Anisotropic mesh adaptation has proved to be a useful tool in numerical solution of partial differential equations (PDEs).
This is especially true when problems arising from science and engineering have distinct anisotropic features. 
The ability to adapt the size, shape, and orientation of mesh elements according to certain quantities of interest can significantly improve the accuracy of the solution and enhance the computational efficiency. 

Criteria for an optimal anisotropic triangular mesh were already given by D'Azevedo~\cite{D'Azev91} and Simpson~\cite{Simpso94} in the early nineties of the last century. 
A number of algorithms for automatic construction of such meshes have since been developed. 

A common approach for generating an anisotropic mesh is based on generation of a quasi-uniform mesh in some metric space. 
A key component of the approach is the determination of an appropriate metric often based on some type of error estimates. 
Unfortunately, classic isotropic error estimates do not suit this purpose well because they generally do not take the directional effect of the error or solution derivatives into consideration.
This explains the recent interest in anisotropic error estimation; for example, see anisotropic interpolation error estimates by Formaggia and Perotto~\cite{ForPer01}, Huang~\cite{Huang05}, and Huang and Sun~\cite{HS01}.
Such error estimates for numerical solution of PDEs can be found, among others, in works by Apel~\cite{Apel99}, Kunert~\cite{Kunert01}, Formaggia and Perotto~\cite{ForPer03}, and Picasso~\cite{Picass03}.

It is worth pointing out that most existing anisotropic error estimates are a priori, requiring information of the exact solution of either the underlying problem or its adjoint, which is typically unavailable in a numerical simulation. 
A widely-used approach of avoiding this difficulty in practical computation is to replace the information by one recovered from the obtained numerical approximation.
A number of recovery techniques can be used for this purpose, such as the gradient recovery technique by Zienkiewicz and Zhu~\cite{ZieZhu92,ZieZhu92a} and  the technique based on the variational formulation by Dolej\v{s}{\'\i}~\cite{Dolejs98}. 
Zhang and Naga~\cite{ZhaNag05} have recently proposed a new algorithm to reconstruct the gradient (which can also be used to reconstruct the Hessian) by fitting a quadratic polynomial to the nodal function values and subsequently differentiating it. 
It has been shown by Zhang and Naga~\cite{ZhaNag05} and by Vallet et al.~\cite{VaMDDG07} that the latter is robust and works best among several recovery techniques.
Generally speaking, recovery methods work well when exact nodal function values are used but may lose some accuracy when applied to finite element approximations on non-uniform meshes.
Nevertheless, the optimality of mesh adaptation based on those recovered approximations can still be proven under suitable conditions, see Vassilevski and Lipnikov~\cite{VasLip99}.
More recently, conditions for asymptotically exact gradient and convergent Hessian recovery from a hierarchical basis error estimator have been given by Ovall~\cite{Ovall07}. 
His result is based on superconvergence results by Bank and Xu~\cite{BanXu03,BanXu03a}, which require that the mesh be uniform or almost uniform.

The objective of this paper is to study the use of a posteriori error estimates in anisotropic mesh adaptation.
Although a posteriori error estimates are frequently used for mesh adaptation, especially for refinement strategies and recently also for construction of
equidistributing meshes for numerical solution of two-point boundary value problems by He and Huang~\cite{HH08a} as well as in connection with the moving finite element method by Lang et al.~\cite{LaCaHR03}, up to now only few methods for their use in anisotropic mesh adaptation have been published. 
For example, Cao et al.~\cite{CaHuRu01} studied two a posteriori error estimation strategies for computing scalar monitor functions for use in adaptive mesh movement; Apel et al.~\cite{ApGrJM04} investigated a number of a posteriori strategies for computing error gradients used for directional refinement; and Agouzal et al.~\cite{AgLiVa08} proposed a new method for computing tensor metrics provided that an edge-based a posteriori error estimate is given.
Moreover, Dobrowolski et al.~\cite{DoGrPf99} have pointed out that error estimation based on solving local error problems can be inaccurate on anisotropic meshes.
This shortcoming of local error estimates can be explained by the fact that they generally do not contain enough directional information of the solution, which is global in nature, and that their accuracy and effectiveness are sensitive to the aspect ratio of elements, which can be large for anisotropic meshes. 
We thus choose to develop our approach based on error estimation by means of globally defined error problem.
To enhance the computational efficiency, we employ an iterative method to obtain a cost-efficient approximation to the solution of the corresponding global linear system.
Numerical results show that a few symmetric Gauß--Seidel iterations are sufficient for this purpose.
This is not surprising since the approximation is used only in mesh generation and it is often unnecessary to compute the mesh to a very high accuracy as for the solution of the underlying differential equation.
Numerical experiments also show that the new approach is comparable in accuracy and efficiency to methods using Hessian recovery.
We also test it with a more challenging example: a heat conduction problem for a thermal battery with large and orthotropic jumps in the material coefficients.\footnote{A Sandia National Laboratories benchmark problem.} 

The outline of the paper is as follows. In \cref{sec:FE}, the new framework of using a posteriori hierarchical error estimates for anisotropic mesh adaptation in finite element approximation is described.
In \cref{sec:MetricTensor}, the optimal metric tensor based on the interpolation error is developed.
Several implementation issues are addressed in \cref{sec:Computation}.
Numerical results obtained with the new approach and with Hessian recovery-based methods are presented in \cref{sec:Examples} for a selection of test examples. 
Numerical results for the heat conduction problem are given in \cref{sec:Real-Life}. 
Finally, \cref{sec:Conclusion} contains conclusions and comments.

\section{Model problem and adaptive finite element approximation}\label{sec:FE}

In this section, we describe a new framework of using a posteriori 
hierarchical error estimates for anisotropic mesh adaptation in finite element approximation.

\subsection{Model problem and finite element approximation}
Consider the boundary value problem of a second-order elliptic differential equation.
Assume that the corresponding variational problem is given by
\begin{align*}
   (P)\qquad \begin{cases}
  \text{Find } u\in V \text{ such that} \\
     a(u,v) = F(v), \quad \forall v \in V
   \end{cases}
\end{align*}
where $V$ is an appropriate Hilbert space of functions over a domain $\Omega \in \R^2$, $a(\cdot,\cdot)$ is a bilinear form defined on $V\times V$, and $F(\cdot)$ is a continuous linear functional on $V$.
The finite element approximation $u_h$ of $u$ is the solution of the corresponding variational problem on a finite dimensional subspace $V_h \subset V$, i.e.,
\begin{align*}
  (P_h) \qquad \begin{cases}
  \text{Find } u_h\in V_h \text{ such that} \\
     a(u_h,v_h) = F(v_h), \quad \forall v_h \in V_h.
   \end{cases}
\end{align*}
If the bilinear form $a(\cdot, \cdot)$ is coercive and continuous on $V$, both variational problems $(P)$ and $(P_h)$ have unique solutions. 
The finite dimensional subspace $V_h$ is often chosen as a space of piecewise polynomials associated with a given mesh, say $\T_h$,  on $\Omega$. The variational problem $(P_h)$ results in a system of $\dim(V_h)$ linear algebraic equations. 

\subsection{Adaptive linear finite element solution}

In this work we consider a linear finite element method, where $V$ is taken as $H^1(\Omega)$ and $V_h$ is the space of continuous, piecewise linear functions over $\T_h$.

Let $\T_h^{(i)}$ ($i=0, 1, \dotso$) be an affine family of simplicial meshes on $\Omega$ and $V_h^{(i)}$ the corresponding space of continuous, piecewise linear functions.
The adaptive solution is the result of an iterative process described as follows. 

\begin{figure}[t]
   \hfill{}
   \subcaptionbox{Surface plot.\label{fig:anisotropic_example_a}}{
      \includegraphics[width=0.35\textwidth,clip]{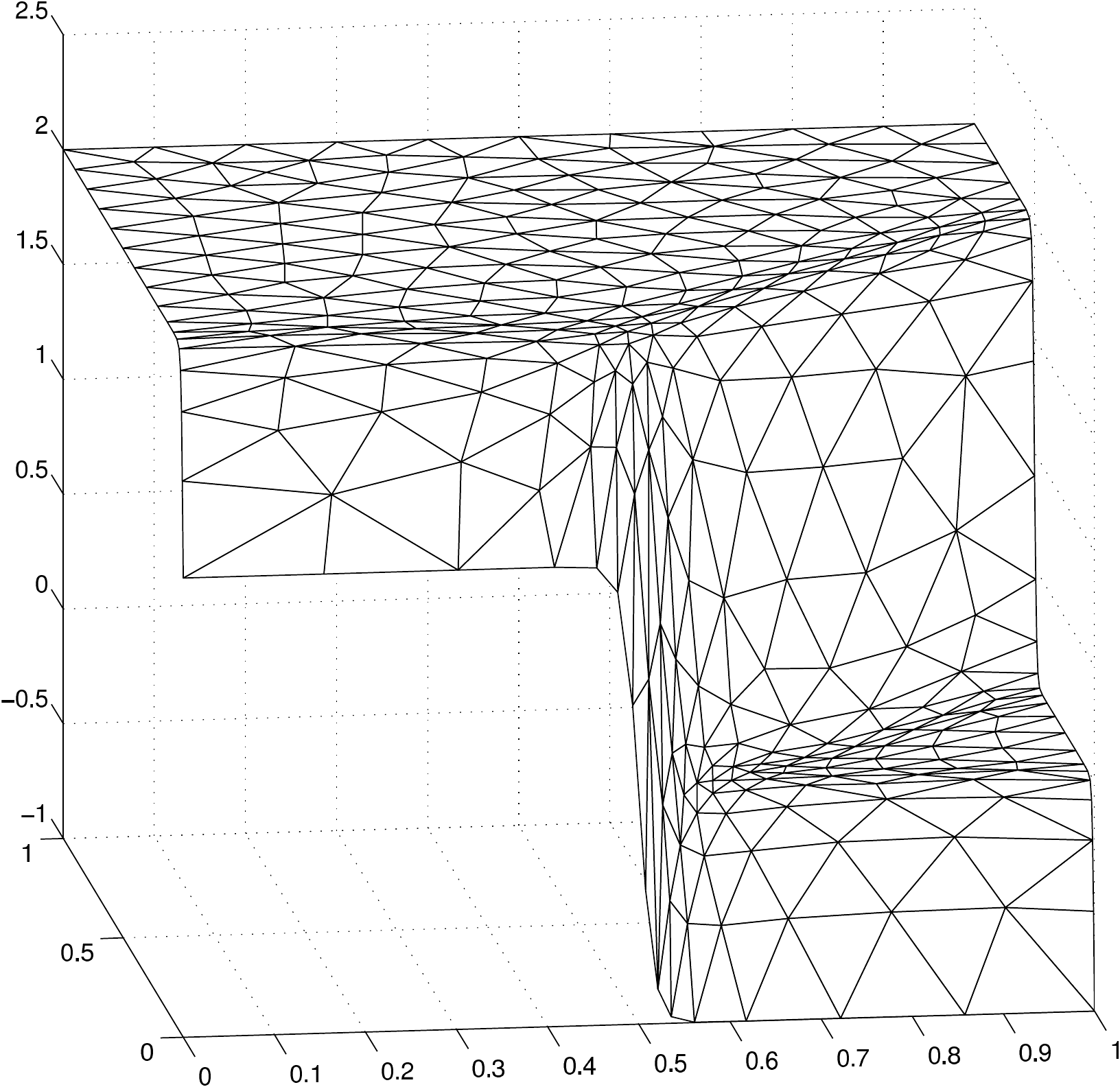}
   }
   \hfill{}
   \subcaptionbox{The corresponding adaptive mesh.\label{fig:anisotropic_example_b}}{
      \includegraphics[width=0.35\textwidth,clip]{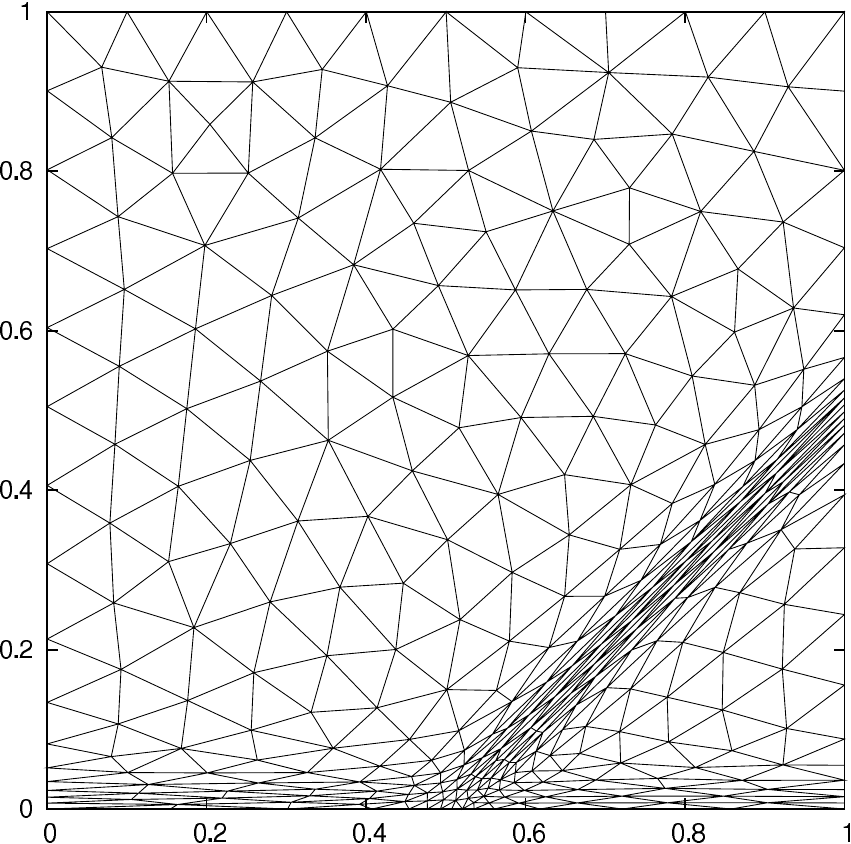}
   }
   \hfill{}
   \caption{An example of anisotropic mesh adaptation for the test 
      function  $u(x,y) = \tanh(60x) - \tanh\left(60(x-y)-30\right)$: 
      Surface plot (\protect\subref{fig:anisotropic_example_a}) of the function 
      on an adaptive mesh
      (\protect\subref{fig:anisotropic_example_b}) obtained with the use
      of the exact Hessian.}\label{fig:anisotropic_example}
\end{figure}

We start with an initial mesh $\T_h^{(0)}$. On every mesh $\T_h^{(i)}$ we solve the variational problem $(P_h)$ with $V_h^{(i)}$ and use the obtained approximation $u_h^{(i)}$ to compute a new adaptive mesh for the next iteration step.
The new mesh $\T_h^{(i+1)}$ is generated as an almost uniform mesh in a metric space with a metric tensor $M_h^{(i)}$ defined in terms of $u_h^{(i)}$.
This yields the sequence
\[ (\T_h^{(0)}, V_h^{(0)}) \rightarrow u_h^{(0)} \rightarrow M_h^{(0)} \rightarrow 
   (\T_h^{(1)}, V_h^{(1)}) \rightarrow u_h^{(1)} \rightarrow M_h^{(1)} \rightarrow \cdots
\]
The process is repeated until a good adaptation is achieved.
An example of such adaptive meshes is shown in \cref{fig:anisotropic_example}.

Typically, the metric tensor $M_h$ depends on the Hessian of the exact solution of the underlying
problem~\cite{ForPer01, Huang05a}, which is often unavailable in practical computation.
The common approach to avoid this difficulty is to recover an approximate Hessian from the computed solution.
We consider here an alternative approach, which uses an a posteriori error estimator for defining and
computing $M_h$.

\subsection{Mesh adaptation based on a posteriori error estimates}
Let $R_h$ be a reconstruction operator applied to the numerical approximation $u_h$. It can be either
a recovery process, a smoothing operator, or an operator connected to an a posteriori error estimate.
We assume that the reconstruction $R_h u_h$ is better than $u_h$ in the sense that
\begin{align}
   \|R_h u_h - u\| \leq \beta \|u_h - u\|,
   \label{eq:sa}
\end{align}
for a given norm $\|\cdot\|$, where $0 \leq \beta < 1$ is a constant.

From the triangle inequality we immediately have
\begin{align} 
   \|u-u_h\| \le\frac{1}{1-\beta} \|R_h u_h-u_h\|. 
   \label{eq:rhuh-uh} 
\end{align}
If the reconstruction $R_h$ has the property
\begin{equation}
I_h R_h v_h = v_h\qquad \forall v_h \in V_h
\label{eq:ihrh}
\end{equation}
for some interpolation operator $I_h$,
we can bound the finite element approximation error by the (explicitly computable) interpolation error of
the reconstructed function $R_h u_h$, viz.,
\begin{align}
   \|u-u_h\| 
     \le \frac{1}{1-\beta} \|R_h u_h-u_h\|
      =  \frac{1}{1-\beta}\|R_h u_h-I_h R_h u_h\|.
   \label{eq:ehrhuh}
\end{align}
Moreover, from the interpolation theory we know that the interpolation error for a given function $v$ can be bounded by a term depending on the triangulation $\T_h$ and derivatives of $v$, i.e.,
\begin{equation}
\| v - I_h v \| \leq C \cdot  \mathcal{E}(\T_h,v),
\label{eq:interperr}
\end{equation}
where $C$ is a constant independent of $\T_h$ and $v$.
Therefore, we can rewrite \cref{eq:ehrhuh} as
\begin{align}
   \|u-u_h\| \le\frac{C}{1-\beta} \; \mathcal{E}(\T_h,R_h u_h) .
   \label{eq:ie_th_rhuh}
\end{align}
In other words, up to a constant, the solution error is bounded by the interpolation error of $R_h u_h$.

\subsection{Hierarchical basis}
One possibility to achieve the property \cref{eq:ihrh} is to use the hierarchical decomposition
of the finite element space. Let 
   \[
   \bar V_h = V_h \oplus W_h,
    \]
where $W_h$ is a hierarchical extension of $V_h$ to $\bar V_h$. Each $\bar v_h \in \bar V_h$ has
a unique representation $\bar v_h = v_h + w_h$ with $v_h \in V_h$ and $w_h \in W_h$. 
If an interpolation operator, $\bar I_h:\bar V_h \mapsto V_h$, can be defined such that
\begin{equation}
   \bar I_h w_h = 0,\quad \forall w_h \in W_h
   \label{eq:bih}
\end{equation}
and if we define $R_h$ through
\begin{equation}
   R_h u_h = u_h + z_h
   \label{eq:rh}
\end{equation}
for some $z_h \in W_h$, then we shall have the property \cref{eq:ihrh} and the estimate \cref{eq:ie_th_rhuh}.
Moreover, 
\begin{align*}
  \|R_h u_h - \bar I_h R_h u_h\| = \|u_h+z_h - u_h\| = \|z_h\| = \|z_h - \bar I_h z_h\| .
\end{align*}
Consequently, we can estimate the finite element approximation error by evaluating
the interpolation error of $z_h$, i.e.,
\begin{align}
  \|u-u_h\| \le 
       \frac{1}{1-\beta} \|z_h - \bar I_h z_h\|
       \leq \frac{C}{1-\beta} \; \mathcal{E}(\T_h,z_h). 
      \label{eq:fe_error_eh}
\end{align}
In the context of a posteriori error estimates, $z_h$ is typically taken as a hierarchical basis error estimator.

\subsection{A posteriori error estimate based on hierarchical basis}\label{ssec:hbee}
The computation of the error estimator is based on a general framework, details on which can be found among others in the work of Bank and Smith~\cite{BanSmi93} or Deuflhard et al.~\cite{DeLeYs89}.
The approach is briefly explained as follows.

Let $u_h \in V_h$ be a linear finite element solution of the variational problem $(P_h)$ and let $\bar V_h = V_h \oplus W_h$, where $W_h$ is the linear span of the edge bubble functions. Obviously, $\bar V_h$ is a subspace of piecewise
quadratic functions. 
Moreover, we can define $\bar I_h$ as the vertex-based, piecewise linear Lagrange interpolation. 
This interpolation satisfies \cref{eq:bih} since the edge bubble functions vanish at vertices.

Let $e_h = u - u_h$ be the error of the finite element solution $u_h$. Then for all $v\in V$ we have
\begin{align}
a(e_h,v) &= F(v) - a(u_h,v).
\end{align}
The error estimate $z_h$ is then defined as the solution of the approximate error problem 
\begin{align*}
   (E_h) \qquad \begin{cases}
   \text{Find } z_h\in W_h \text{ such that} \\
   a(z_h, w_h) = F(w_h) - a(u_h,w_h) \quad \forall w_h \in W_h.
   \end{cases}
\end{align*}
The estimate $z_h$ can be viewed as a projection of the true error onto the subspace $W_h$. 
Note that this definition of the error estimate is global and its solution can be costly.
Several solution methods will be discussed in \cref{sec:Computation}. 

Once $z_h$ is determined, the reconstruction $R_h u_h$ is derived from \cref{eq:rh}.
Then, if assumption \cref{eq:sa} holds, the finite element approximation error can be controlled by minimizing the interpolation error of $z_h$, i.e., the right-hand side in \cref{eq:fe_error_eh}.
In this paper, we construct optimal metric tensors with respect to interpolation
error estimates $\mathcal{E}(\T_h,z_h)$ for the $L^2$ norm.
We assume that the reconstruction $R_h u_h = u_h + z_h$, where $z_h$ is computed from ($E_h$), gives a better approximation to $u$ than $u_h$, i.e., $\beta<1$ in \cref{eq:sa}.

\section{Metric tensor based on linear interpolation error estimate}\label{sec:MetricTensor}

\subsection{Equidistribution and alignment}
Let $\Omega$ be a polyhedral domain in $\R^d$ and let $\T_h$ be a simplicial triangulation on $\Omega$.
For every element $K\in \T_h$, there exists an affine invertible mapping $F_K \colon  \hat K \to K$
such that $K = F_K(\hat K)$, where $\hat K$ is the reference element. We assume that $\hat K$ has
been chosen to be equilateral and have a unitary volume. We denote the Jacobian matrix of $F_K$ by $F_K'$ and the number of elements in $\T_h$ by $N$.

As mentioned before, we consider an adaptive anisotropic mesh as a uniform mesh in the metric specified
by a metric tensor $M$. 
Such a mesh is referred hereafter to as an {\em $M$-uniform mesh}.
It can be characterized by shape-orientation and size requirements on mesh elements;
see~\cite{Huang06}. 

\emph{Alignment condition (i.e., shape-orientation requirement).} The elements of an $M$-uniform
mesh $\T_h$ are equilateral in the metric specified by $M$. This can be expressed as
\begin{align}
   \frac{1}{d} \tr \left({(F_K')}^T M_K F_K' \right) 
= \det  {\left({(F_K')}^T M_K F_K' \right)}^{\frac{1}{d}}, \quad \forall K \in \T_h
   \label{eq:alignment}
\end{align}
where $M_K$ is the average of $M$ on element $K$,i.e.,
\[
M_K = \frac{1}{|K|}\int_K M(x) d x.
\]
The left-hand side term of equality \cref{eq:alignment} is equal to the arithmetic-mean of
the eigenvalues of matrix ${(F_K')}^T M_K F_K'$ while the right-hand side term is equal to their
geometric-mean. The arithmetic-mean geometric-mean inequality implies that
\cref{eq:alignment} holds if and only if the eigenvalues of matrix ${(F_K')}^T M_K F_K'$
are all equal. Element $K$ is equilateral in the metric $M_K$
when it satisfies \cref{eq:alignment}.

\emph{Equidistribution condition (i.e., size requirement).} The elements of an $M$-uniform mesh have
an equal volume in the metric $M$, i.e.,
\begin{align}
   |K| \sqrt{\det(M_K)} = \frac{\sigma_h}{N}, \quad \forall K \in \T_h
   \label{eq:equidistribution}
\end{align}
where
\[
\sigma_h = \sum\limits_{K\in\T_h}|K|\sqrt{\det(M_{K})}.
\]
Note that the left-hand side of \cref{eq:equidistribution} is equal to the volume of element $K$
in metric $M_K$, i.e.,
\[
\int_K \sqrt{\det(M_K)} d x = |K| \sqrt{\det(M_K)} .
\]

\subsection{Anisotropic interpolation error bound for piecewise quadratic functions}
Elementwise anisotropic interpolation error estimates are developed in~\cite{ForPer01,ForPer03,HS01}. 
Here, we follow the theory in~\cite{HS01}.
Consider the piecewise linear Lagrange interpolation ($k = 1$) of a piecewise quadratic function $v$ on
an arbitrary mesh $\T_h$.
The elementwise interpolation error measured in the $L^q$ norm ($q\ge 1$) is given by
\begin{align*}
   \|v - I_h v \|^q_{L^q(K)} 
&\le C |K| {\left(\tr\left( {(F_K')}^T |H_K| F_K'\right)\right)}^q,
\end{align*}
where $H_K$ is the Hessian of $v$ on the element $K$, $|H_K| = \sqrt{H_K^T H_K}$, $C$
is a constant independent of $\T_h$ and $v$, and $\mbox{tr}(\cdot)$ denotes the trace of a matrix.
Note that $H_K$ is constant on $K$ since by assumption $v$ is quadratic on the element.
Summing over all elements of $\T_h$ provides an upper bound for the global interpolation error 
\begin{equation}
\| v - I_h v \|^q_{L^q(\Omega)} \le C
   \sum_{K \in \T_h} |K| {\left(\tr\left({(F_K')}^T |H_K| F_K'\right)\right)}^q .
\label{eq:interperrbound}
\end{equation}
One may notice that we have used $L^q$ norm for the error. As we shall see later
(cf.~\cref{eq:interperrbound2+1}), an optimal global error bound in this norm can be obtained
for the non-regularized case. In principle, the same procedure also works for other norms
or semi-norms particularly the $H^1$ semi-norm. However, it is unclear that the interpolation
error bounds obtained in~\cite{HS01} for other norms will lead to an optimal global bound
for $M$-uniform meshes.

From this, we can set $\mathcal{E}(\T_h, v)$ in \cref{eq:interperr} to
\begin{align}
   \mathcal{E}(\T_h, v) =
   \sum_{K \in \T_h} |K| {\left(\tr\left({(F_K')}^T \lvert H_K \rvert F_K'\right)\right)}^q .
   \label{eq:global_ie}
\end{align}
It has a lower bound as
\begin{align}
\mathcal{E}(\T_h,v) & = \sum_{K \in \T_h} |K| 
{\left(\tr\left({(F_K')}^T |H_K| F_K'\right)\right)}^q \notag \\
   & \ge d^q\sum_{K \in \T_h} |K| 
   {\left(\det\left({(F_K')}^T |H_K| F_K'\right)\right)}^{\frac{q}{d}}
      \label{eq:eh-first-ineq}\\
      &= d^q\sum_{K \in \T_h} |K|^{\frac{d+2q}{d}} {\det(|H_K|)}^{\frac{q}{d}} 
    \notag \\
    & = d^q\sum_{K \in \T_h} {\left( |K|\; {\det(|H_K|)}^{\frac{q}{d+2q}}\right)}^{\frac{d+2q}{d}} \notag\\
    & \geq d^q N^{-\frac{2q}{d}} {\left( 
    \sum_{K \in \T_h} |K|\; {\det(|H_K|)}^{\frac{q}{d+2q}} \right)}^{\frac{d+2q}{d}},
      \label{eq:eh-second-ineq}
\end{align}
where we have used the arithmetic-mean geometric-mean inequality in \cref{eq:eh-first-ineq}
(recalling the trace and determinant of a matrix are equal to the sum and product of its eigenvalues,
respectively) and H\"older's inequality in \cref{eq:eh-second-ineq}.
If $\max_{K\in\T_h} \diam(K) \rightarrow 0$, where $\diam(K)$ denotes the diameter of $K$,
we see that the asymptotic lower bound on $\mathcal{E}(\T_h,v)$ is
\begin{align}
   d^q N^{-\frac{2q}{d}} {\left( \int_\Omega {\det(|H|)}^{\frac{q}{d+2q}} \; dx \right)}^{\frac{d+2q}{d}} ,
   \label{eq:eh-lowest-bound}
\end{align}
which is invariant for all meshes of the same number of elements $N$.
Thus, a mesh on which $\mathcal{E}(\T_h, v)$ attains a lower bound \cref{eq:eh-second-ineq} can be considered to be an asymptotically optimal mesh.

\subsection{Optimal metric}
The optimal metric $M$ is defined such that the interpolation error bound $\mathcal{E}(\T_h, v)$ defined in \cref{eq:global_ie} attains its lower bound \cref{eq:eh-second-ineq} on $M$-uniform meshes of $N$ elements associated with $M$.

We first notice that equality in \cref{eq:eh-first-ineq} holds if the $M$-uniform mesh satisfies
\begin{align*}
   \frac{1}{d} \tr\left({(F_K')}^T |H_K| F_K'\right)
   = {\det\left({(F_K')}^T |H_K| F_K'\right)}^{\frac{1}{d}}, \quad \forall K \in \T_h .
\end{align*}
Comparing this with the alignment condition \cref{eq:alignment}, a property satisfied by the $M$-uniform mesh,
suggests that $M$ be defined as
\begin{align*}
   M_K = \theta_K |H_K|
\end{align*}
with some scalar function $\theta_K$.

Next we notice that equality in \cref{eq:eh-second-ineq} holds if the mesh satisfies
\begin{align*}
   |K|\; {\det(|H_K|)}^{\frac{q}{d+2q}} 
   = \frac{1}{N}  \sum_{K \in \T_h} |K| \; {\det(|H_{K}|)}^{\frac{q}{d+2q}} ,
         \quad \forall K \in \T_h .
\end{align*}
Comparing this to the equidistribution condition \cref{eq:equidistribution}, another property satisfied
by the $M$-uniform mesh, yields
\begin{align*} 
   \sqrt{\det(M_K)} = {\det\left(|H_K|\right)}^{\frac{q}{d+2q}} .
\end{align*}
This condition can be used for determining $\theta_K$. Thus, we obtain the optimal metric tensor as 
\begin{align} 
   M_K = {\det\left(|H_K|\right)}^{-\frac{1}{d+2q}} |H_K|, \quad \forall K \in \T_h.
   \label{eq:optimal_metric}
\end{align}
The interpolation error bound \cref{eq:global_ie} attains its lower bound \cref{eq:eh-second-ineq} on any $M$-uniform mesh associated with this metric tensor.
From \cref{eq:interperrbound} we obtain
\begin{align}
   \| v - I_h v \|_{L^q(\Omega)} 
      & \le C N^{-\frac{2}{d}} {\left(
      \sum_{K \in \T_h} |K| \; {\det(|H_K|)}^{\frac{q}{d+2q}}
   \right)}^{\frac{d+2q}{d q}}
   \label{eq:interperrbound2} \\
    & \sim  C N^{-\frac{2}{d}} {\left(\int_\Omega {\det(|H|)}^{\frac{q}{d+2q}} \; dx \right)}^{\frac{d+2q}{d q}} 
       \notag \\
    & = C N^{-\frac{2}{d}} \left\|\sqrt[d]{\det(|H|)}\right\|_{L^{\frac{d q}{d+2q}}(\Omega)}
       \label{eq:interperrbound2+1}
\end{align}
for any $M$-uniform mesh associated with the metric tensor \cref{eq:optimal_metric}.
Bound \cref{eq:interperrbound2+1} has been obtained in~\cite{HS01} for $q=2$ and obtained and shown to be optimal in~\cite{ChSuXu07} for general $q\ge 1$.

The metric tensor defined by \cref{eq:optimal_metric} is not necessarily positive definite since both $|H_K|$ and $\det(|H_K|)$ can vanish locally. To avoid this difficulty, the error bound is regularized with a positive parameter $\alpha_h$, i.e.,
\begin{align}
 \| v - I_h v \|^q_{L^q(\Omega)} 
 & \le C   \sum_{K \in \T_h} |K|
 {\left(\frac{1}{d}\tr\left({(F_K')}^T \left [\alpha_h I + |H_K| \right ] F_K'\right)\right)}^q 
\nonumber \\
& = C \alpha_h^q    \sum_{K \in \T_h} |K| {\left(\frac{1}{d}\tr\left({(F_K')}^T \left [I + \frac{1}{\alpha_h}
|H_K| \right ] F_K'\right)\right)}^q .
\label{eq:interperrbound3+1}
\end{align}
Using the same procedure as above, by minimizing the above (regularized) error bound we obtain the optimal metric tensor as
\begin{align} 
   M_K = {\det\left(I+\frac{1}{\alpha_h} |H_K|\right)}^{-\frac{1}{d+2q}} \left (I+\frac{1}{\alpha_h} |H_K|\right )
   , \quad \forall K \in \T_h .
   \label{eq:optimal_metric2}
\end{align}
The regularization parameter plays a role of controlling the intensity of mesh adaptation.
Indeed, as $\alpha_h \to \infty$, $M_K \to I$ and a uniform mesh results.
On the other hand, as $\alpha_h \to 0$, the mesh adaptation is increasingly reliant on $|H_K|$.
To balance between these situations, we follow~\cite{HS01} and define $\alpha_h$
through the algebraic equation
\[ \sum_{K\in \T_h} \sqrt{\det(M_K)}\; |K| = 2^{\max\left\{1,\frac{dq}{d+2q}\right\} } \lvert \Omega \rvert , \]
or equivalently
\begin{equation}
   \sum_{K \in \T_h} 
   {\det\left(I+\frac{1}{\alpha_h} |H_K|\right)}^{\frac{q}{d+2q}} \; |K| 
      = 2^{\max\left\{1,\frac{dq}{d+2q}\right\} }  \lvert \Omega \rvert ,
\label{eq:alpha}
\end{equation}
where the factor $2^{\max\left\{1,\frac{dq}{d+2q}\right\} }$ has been used so that
lower and upper bounds can be obtained for $\alpha_h$; see \cref{alpha-bound} and its derivation below.
With this definition, about half of the mesh elements are concentrated in regions where
$\det(M_K)$ is large~\cite{HS01}. Moreover, $M_K$ is invariant under a scaling transformation
of $v$.

Equation \cref{eq:alpha} has a unique solution since its left-hand side is monotonically decreasing
with $\alpha_h$ increasing (assuming that $|H_K|$ is not all zero for all elements of $\mathcal{T}_h$),
and tends to $+\infty$ (which is greater than the right-hand side)
as $\alpha_h \to 0$ and $|\Omega|$  (which is less than the right-hand side) as $\alpha_h \to \infty$.
Moreover, it can be solved using a simple
iteration scheme such as the bisection method.
Furthermore, lower and upper bounds on $\alpha_h$ can be obtained,
\begin{align}
   & {\left[{\left(2^{\max\left\{2,\frac{dq}{d+2q}+1\right\} -\frac{q}{d+2q}} - 1\right)}^{-1}
   |\Omega|^{-1} \sum_{K \in \T_h} {\det\left( |H_K|\right)}^{\frac{q}{d+2q}} \; |K|
\right]}^{\frac{d+2q}{d q}} 
\nonumber \\
&\qquad \mbox{ } \qquad \le \alpha_h \le
{\left[\frac{1}{ |\Omega|} \sum_{K \in \T_h} \|H_K\|^{\frac{d q}{d+2q}} \; |K|
\right]}^{\frac{d+2q}{d q}}.
\label{alpha-bound}
\end{align}
Indeed, from \cref{eq:alpha} we have
\begin{align}
2^{\max\{1,\frac{dq}{d+2q}\} } |\Omega|
& = \sum_{K \in \T_h} {\det\left(I+\frac{1}{\alpha_h} |H_K|\right)}^{\frac{q}{d+2q}} \; |K| 
\notag \\
& \le \sum_{K \in \T_h} {\left\| I+\frac{1}{\alpha_h} |H_K|\right\|}^{\frac{dq}{d+2q}} \; |K| 
\notag \\
& \le \sum_{K \in \T_h} {\left( 1 +\alpha_h^{-1} \| H_K\| \right)}^{\frac{dq}{d+2q}} \; |K| 
\notag \\
& \le 2^{\max\{0,\frac{dq}{d+2q}-1\}} \sum_{K \in \T_h} \left( 1 +\alpha_h^{-\frac{dq}{d+2q}}
 \| H_K\|^{\frac{dq}{d+2q}}\right) \; |K| 
\notag \\
& =  2^{\max\{0,\frac{dq}{d+2q}-1\}} 
\left(|\Omega| + \alpha_h^{-\frac{dq}{d+2q}} \sum_{K \in \T_h} \| H_K\|^{\frac{dq}{d+2q}} \; |K| \right),
\notag
\end{align}
which leads to the right inequality of \cref{alpha-bound}. On the other hand,
\begin{align}
2^{\max\{1,\frac{dq}{d+2q}\} } |\Omega|
& \ge \sum_{K \in \T_h} {\left(1+\alpha_h^{-d} \det\left (|H_K|\right)\right)}^{\frac{q}{d+2q}} \; |K| 
\notag \\
& \ge 2^{\frac{q}{d+2q}-1} \sum_{K \in \T_h} \left(1+\alpha_h^{-\frac{dq}{d+2q}}
   {\det\left (|H_K|\right)}^{\frac{q}{d+2q}} \right ) \; |K| 
\notag \\
& = 2^{\frac{q}{d+2q}-1}  \left(|\Omega|+\alpha_h^{-\frac{dq}{d+2q}}
   \sum_{K \in \T_h} {\det\left (|H_K|\right)}^{\frac{q}{d+2q}}  \; |K| \right),
\notag
\end{align}
which gives the left inequality of \cref{alpha-bound}.

The interpolation error bound for a corresponding $M$-uniform mesh can be obtained as follows.
From \cref{eq:interperrbound3+1} and using the equidistribution and alignment conditions we have
\begin{align}
& \| v - I_h v \|^q_{L^q(\Omega)} 
\nonumber \\
& \le C \alpha_h^q  \sum_{K \in \T_h} |K| \; {\det \left(I + \frac{1}{\alpha_h} |H_K| \right)}^{\frac{q}{d+2q}}
{\left(\frac{1}{d}\tr\left({(F_K')}^T M_K F_K'\right)\right)}^q
\nonumber \\
& = C \alpha_h^q  \sum_{K \in \T_h} |K| \; {\det(M_K)}^{\frac{1}{2}}
\; {\det\left({(F_K')}^T M_K F_K'\right)}^{\frac{q}{d}}
\nonumber \\
& = C \alpha_h^q  \sum_{K \in \T_h} |K| \; {\det(M_K)}^{\frac{1}{2}} \;
{\left( |K| \; {\det(M_K)}^{\frac{1}{2}}\right)}^{\frac{2 q}{d}}
\nonumber \\
& = C \alpha_h^q N^{-\frac{2q}{d}} \sigma_h^{\frac{d+2q}{d}} .
\notag
\end{align}
For $\alpha_h$ defined in \cref{eq:alpha}, $\sigma_h = 2^{\max\left\{1,\frac{dq}{d+2q}\right\} } |\Omega|$. Combining this with
\cref{eq:interperrbound3+1} we obtain
\begin{align}
\| v - I_h v \|_{L^q(\Omega)} & \le C N^{-\frac{2}{d}} \alpha_h .
      \label{eq:interperrbound3}
\end{align}

In our computation we use the mesh generation software \emph{bamg} (\emph{bidimensional anisotropic mesh generator} developed by F. Hecht~\cite{bamg}) to generate new adaptive meshes for a given metric tensor $M$.
Note that \emph{bamg} requires that the metric tensor be further normalized such that all elements have a unitary volume in the metric.
Thus, in actual computation we use a normalized metric tensor
\begin{align} 
   \mathcal{M}_K = {\left(\frac{\sigma_h}{N}\right)}^{-\frac{2}{d}}
   \det {\left(I+\frac{1}{\alpha_h} |H_K|\right)}^{-\frac{1}{d+2q}}
   \left(I+\frac{1}{\alpha_h} |H_K|\right),
   \label{eq:optimal_metric2+1}
\end{align}
where $N$ is the desired number of mesh elements and
\[
   \sigma_h = \sum_{K\in \mathcal{T}_h} |K| ~ {\det(M_K)}^{\frac{1}{2}}
= \sum_{K\in \mathcal{T}_h} |K| \det {\left(I+\frac{1}{\alpha_h} |H_K|\right)}^{\frac{q}{d+2q}}  .
\]
It is remarked that the metric tensor can also be normalized using a prescribed error level; see~\cite{Huang05a}.

\section{Computation of the metric tensor and anisotropic meshes}\label{sec:Computation}
We discuss here some implementation issues for two-dimensional problems.

The computation typically starts with a coarse regular Delaunay mesh of the domain and a desired number of mesh elements, $N$. 
For a given triangular mesh $\T_h^{(i)}$ at step $i$, we compute the numerical approximation $u_h^{(i)}$ with a standard linear finite element method. 
Based on $u_h^{(i)}$ and $\T_h^{(i)}$, we then compute $z_h^{(i)}$ as an approximation to the solution of the approximate error problem $(E_h)$.
Once $z_h^{(i)}$ has been obtained, it is straightforward to compute its elementwise Hessian and define the new metric tensor $M^{(i)}$ according to \cref{eq:optimal_metric2}, 
\begin{align*}
   M_K^{(i)} = 
   \det{\left(I + \frac{1}{\alpha_h^{(i)}} |H_K(z_{h}^{(i)})|\right)}^{-\frac{1}{6}}
      \left(I + \frac{1}{\alpha_h^{(i)}} |H_K(z_{h}^{(i)})|\right) ,
\end{align*}
where the error is measured in the $L^2$-norm, i.e., $q=2$.
A new mesh is generated with \emph{bamg} according to the metric tensor $\mathcal{M}^{(i)}= {\left(\sigma_h^{(i)} / N\right)}^{-1} M^{(i)}$. The process is repeated until a good adaptation (see discussion below) is achieved.

\subsection{Mesh quality measure}
In order to characterize the mesh adaptation quality and to define an appropriate stopping criterion for the mesh adaptation process, we introduce the alignment and equidistribution quality measures~\cite{Huang05}
\begin{align*}
Q_{ali}^{(i)}(K) &= {\left[\frac{ \tr \left({(F_K')}^T M_K^{(i)} F_K' \right)}
{d \det  {\left({(F_K')}^T M_K^{(i)} F_K' \right)}^{\frac{1}{d}}} \right]}^{\frac{d}{2(d-1)}}
\end{align*}
and
\begin{align*}
   Q_{eq}^{(i)}(K)& =  \frac{N^{(i)} |K|.
   \sqrt{\det(M_K^{(i)})}} {\sigma_h^{(i)}}, 
\end{align*}
which characterize how closely the mesh satisfies the alignment and equidistribution conditions \cref{eq:alignment} and \cref{eq:equidistribution}, respectively.

Using $M^{(i)}$, $Q_{ali}^{(i)}$, and $Q_{eq}^{(i)}$, the estimate \cref{eq:interperrbound3+1} can be reformulated as
\begin{align*}
   \| v - I_h v \|_{L^q(\Omega)}
   & \leq   C \alpha_h^{(i)} {\left( \sum_{K \in \T_h} |K| 
      {\left(\frac{1}{d}\tr\left({(F_K')}^T \left [I + \frac{1}{\alpha_h} |H_K| \right ] 
         F_K'\right)\right)}^q \right)}^{\frac{1}{q}}\\
   &= C N^{-\frac{2}{d}} \alpha_h^{(i)} 
         {\left(\sigma_h^{(i)}\right)}^{\frac{d+2q}{qd}}
         {\left( \frac{1}{\sigma_h^{(i)}}
         \sum_{K \in \T_h} |K| ~  \sqrt{\det(M_K^{(i)})} 
         {\left(Q_{ali}^{(i)}(K)\right)}^q
         {\left(Q_{eq}^{(i)}(K)\right)}^{\frac{2q}{d}}
      \right)}^{\frac{1}{q}}\\
      & = C  \, N^{-\frac{2}{d}}  \alpha_h^{(i)} {\left(\sigma_h^{(i)}\right)}^{\frac{d+2q}{qd}} Q_{mesh}^{(i)},
\end{align*}
where
\begin{align*}
   Q_{mesh}^{(i)} \equiv {\left[ \frac{1}{\sigma_h^{(i)}} 
      \sum_{K \in \T_h} |K| 
      \, \sqrt{\det(M_K^{(i)})} \, 
      {\left(Q_{ali}^{(i)}(K)\right)}^q
      \, {\left(Q_{eq}^{(i)}(K)\right)}^{\frac{2q}{d}}
      \right]}^{\frac{1}{q}}
\end{align*}
is the overall mesh quality measure and takes into account both the shape and the size of elements. 
 Since $Q_{ali}$ and $Q_{eq}$ appear in $Q_{mesh}$ as a product, their effects are not independent
 but compensate for each other.
 As a consequence, the mesh can have a good overall quality when small elements are shaped worse than large elements or well-aligned elements are worse shaped than worse aligned elements.
Note that $Q_{ali}$, $Q_{eq}$, $Q_{mesh}$ $\geq 1$;
and $Q_{ali} = Q_{eq} = Q_{mesh} = 1$ if and only if the underlying mesh is $M$-uniform (cf.~\cref{eq:interperrbound3}).

In the following numerical tests, the mesh adaptation process is stopped when
\[ Q_{mesh}^{(i)} \leq 1 + \varepsilon_{mesh}, \]
where $\varepsilon_{mesh}$ is a tolerances chosen as $\varepsilon_{mesh}=0.1$ in our computation.

\subsection{Computation of the error estimator}
A key component of the procedure is to find the solution $z_h$ of problem $(E_h)$.
Note that $(E_h)$ is a global problem and finding its exact solution can be as costly as for computing
a quadratic finite element approximation to the original PDE problem. Three approaches are considered
here for solving or approximating $(E_h)$.

\emph{Edge-based error estimator.} The expense of the error estimation can be significantly reduced, if the bilinear form $a$ in $(E_h)$ is replaced by an approximation $\tilde a$ that allows a more efficient solution of the resulting linear system. 
A very efficient approach in two dimensions is to reduce the original problem to a series of local error problems which are defined over two elements sharing a common edge and can be solved efficiently.
The approach is equivalent to the application of one Jacobi's iteration (starting from zero) to the linear system resulting from the global error problem, i.e to the replacement of the stiffness matrix resulting from $(E_h)$ by its diagonal.
This approach has been successfully used in finite element computations~\cite{DeLeYs89, Lang01, LaCaHR03}.
Moreover, it has been shown~\cite{DeLeYs89} that such an error estimator is spectrally nearly equivalent to the original one under suitable conditions.

Despite its success in isotropic mesh adaptation, the approach does not seem to work well for anisotropic mesh adaptation. 
This may be explained by the fact that estimators based on local error problems generally depend on the aspect ratio of elements and can become inaccurate when the aspect ratio is large, a case that is often true for anisotropic meshes. 
Moreover, such estimators may not contain enough directional information of the solution which is global in nature and essential to the success of anisotropic mesh adaptation.

\emph{Node-based error estimator.} This approach is similar to the edge-based error estimator, with the error estimator being obtained by solving a series of local error problems defined on node patches with homogeneous Dirichlet boundary conditions.

\emph{Inexact solution of the full error problem.} In this approach the full error problem is kept but only an approximation to its exact solution is sought and used for the computation of the metric tensor.
In our experiments, a few symmetric Gauß--Seidel iterations are employed to obtain such an approximation. 
In the following computation, Gauß--Seidel iterations are repeated until the relative difference of the old and the new approximations is under a given tolerance \emph{GS-RTOL}.

It is noted that globally defined error estimators have the advantages that they are often independent of element aspect ratio and contain more directional information of the solution.
Moreover, it is known~\cite{DoGrPf99} that the full hierarchical basis error estimator is efficient and reliable for anisotropic meshes.

Numerical comparison among these approaches is given in the next section.

\section{Numerical examples}\label{sec:Examples}
In this section, we present some numerical results for a selection of two-dimensional problems with an anisotropic behaviour. 
We first compare different approaches in solving the error problem $(E_h)$ and then the new method with some common Hessian recovery methods. 
At the end of the section, we give further examples to demonstrate the ability of the method to generate appropriate anisotropic meshes.

Convergence is illustrated by plotting the finite element solution error against the number of elements. 
We use the $L^2$-norm for the error because the monitor function $M_K$ is optimized for this norm.
For the inexact solution of the full error problem, $\text{\emph{GS-RTOL}} = 0.01$ is chosen as a relative tolerance for the iterative Gauss--Seidel approximation.
\subsection{A first example}\label{ex:tanh}
Consider the boundary value problem
\begin{equation}
\begin{cases}
  -\Delta  u = f  & \text{in } \Omega, \\
           u = g  & \text{on } \partial \Omega
\end{cases}
\label{eq:bvp}
\end{equation}
with $\Omega = (0,1)\times(0,1)$. The right-hand side $f$ and the Dirichlet boundary conditions are chosen such that the exact solution is given by
\begin{align*}
   u(x,y) &= \tanh(60x) - \tanh\left(60(x-y)-30\right).
\end{align*}
The solution exhibits a strong anisotropic behaviour and describes the interaction between a boundary layer along the $x$-axis and a shock wave along the line $y = x-0.5$.
A solution plot is given in \cref{fig:anisotropic_example_a}.

\textbf{Reduced vs.\ full error estimators.}
As mentioned in \cref{sec:Computation}, on anisotropic meshes, there can be a significant difference in accuracy between estimators obtained by solving localized error problems and those obtained by means of a globally defined error problem.
In our first test, we investigate the influence of the three error estimators described in the previous section on mesh adaptivity.

\begin{figure}[t] \centering 
   \includegraphics[width=0.65\textwidth,clip]{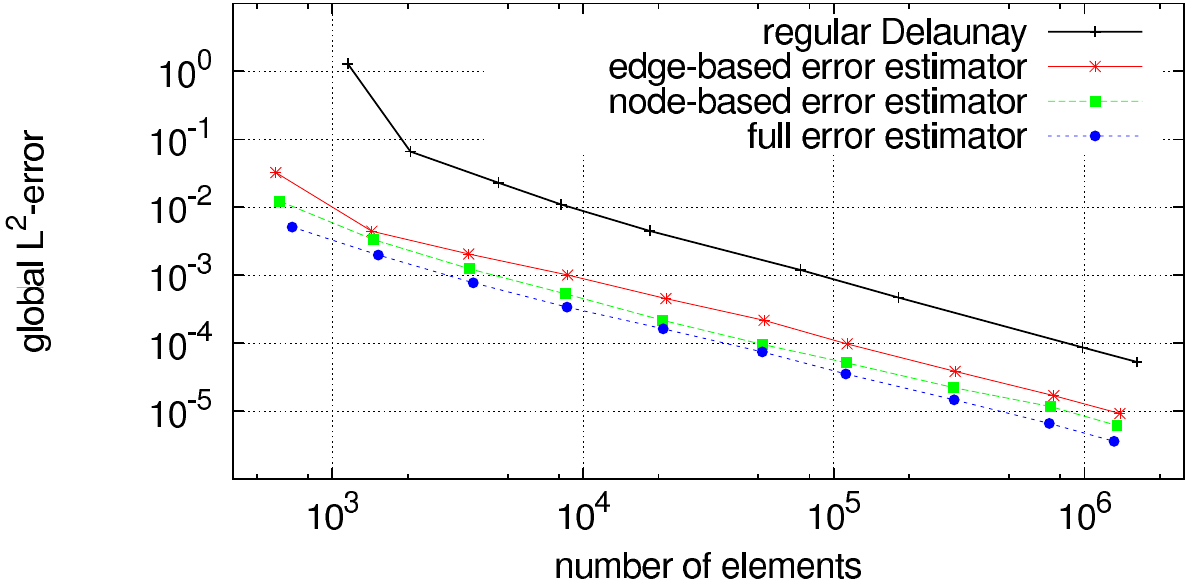}
   \caption{Example~\ref{ex:tanh}: a comparison of the error for adaptive
      finite element solutions obtained with mesh adaptation controlled by
      the reduced and full error estimators.}\label{fig:tanh_eh_comparison}
\end{figure}

\begin{figure}[p]\centering{}%
   \subcaptionbox{Edge-based error estimator: 
      722 vertices and \num{1345} triangles, 
   $\|e\|_{L^2} = \num{5.0e-3}$, maximum aspect ratio $12.8$.\label{fig:tanh_eh_meshes_diagonal}}[1.0\linewidth]{
      \includegraphics[height=0.28\textheight,clip]{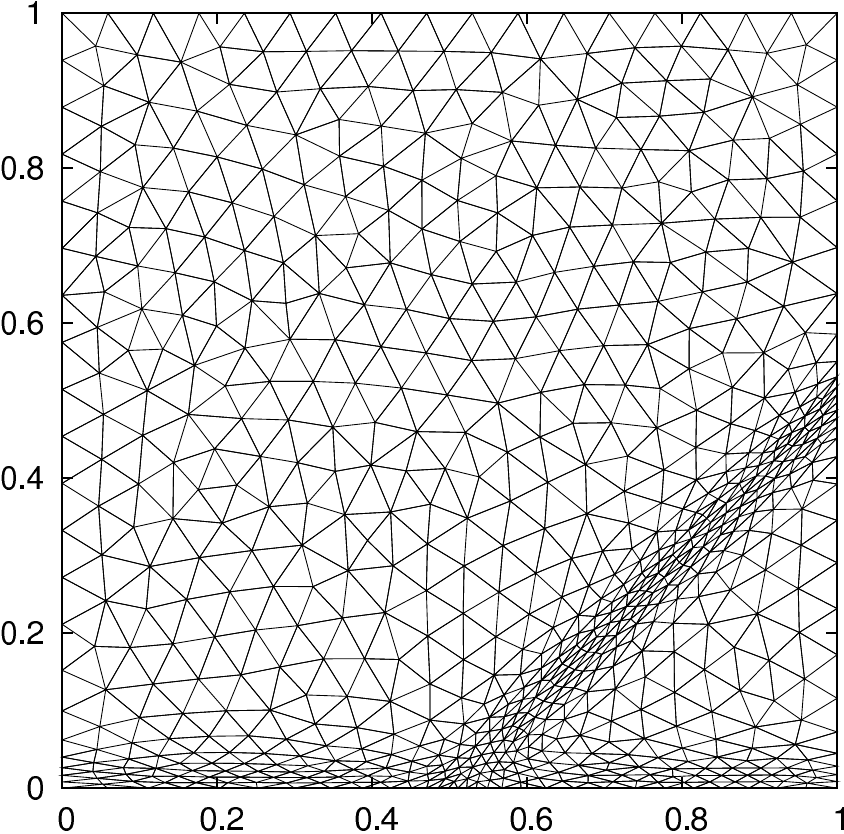}
      \qquad\qquad
      \includegraphics[height=0.28\textheight,clip]{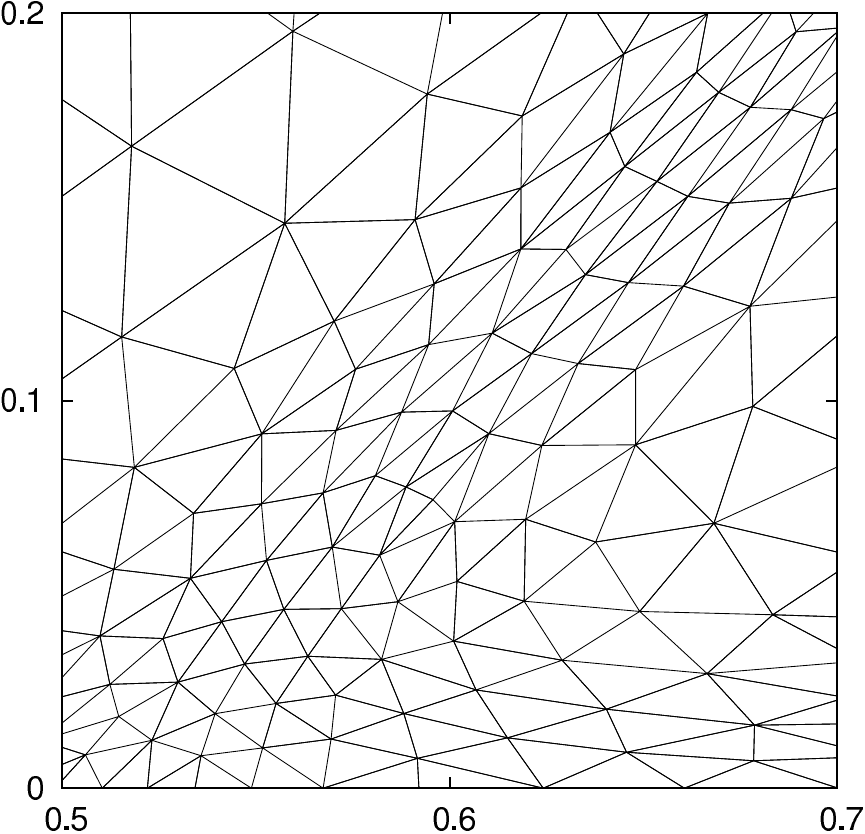}
   }%
   \\[2ex]
   \subcaptionbox{Node-based error estimator: 
      720 vertices and \num{1339} triangles, 
      $\|e\|_{L^2} = \num{3.7e-3}$, maximum aspect ratio $14.3$.\label{fig:tanh_eh_meshes_patch}}[1.0\linewidth]{
      \includegraphics[height=0.28\textheight,clip]{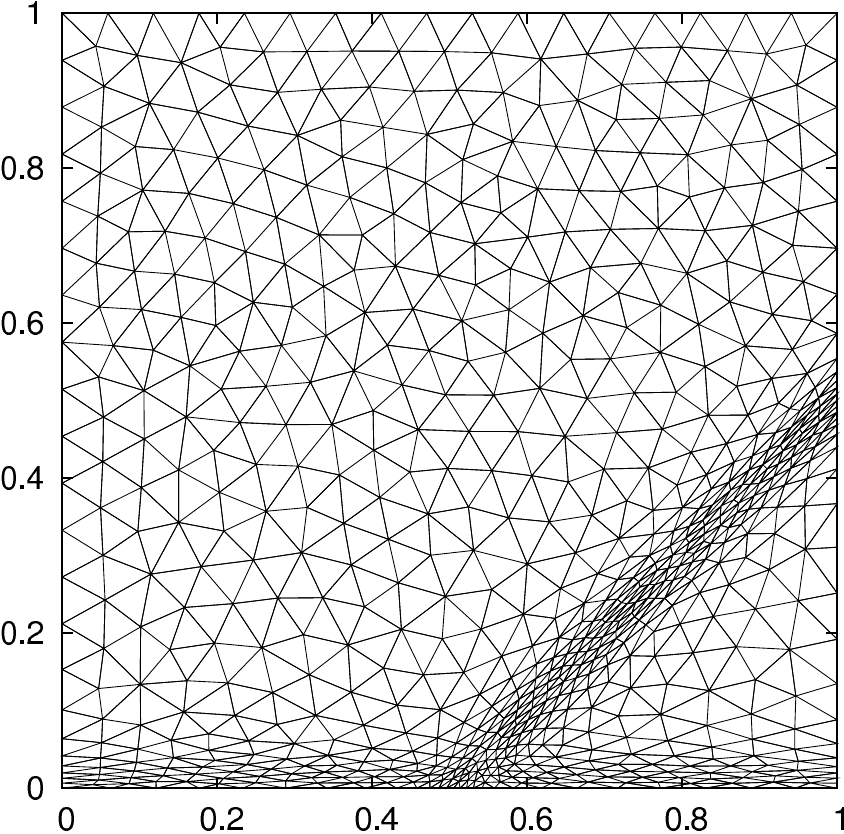}
      \qquad\qquad
      \includegraphics[height=0.28\textheight,clip]{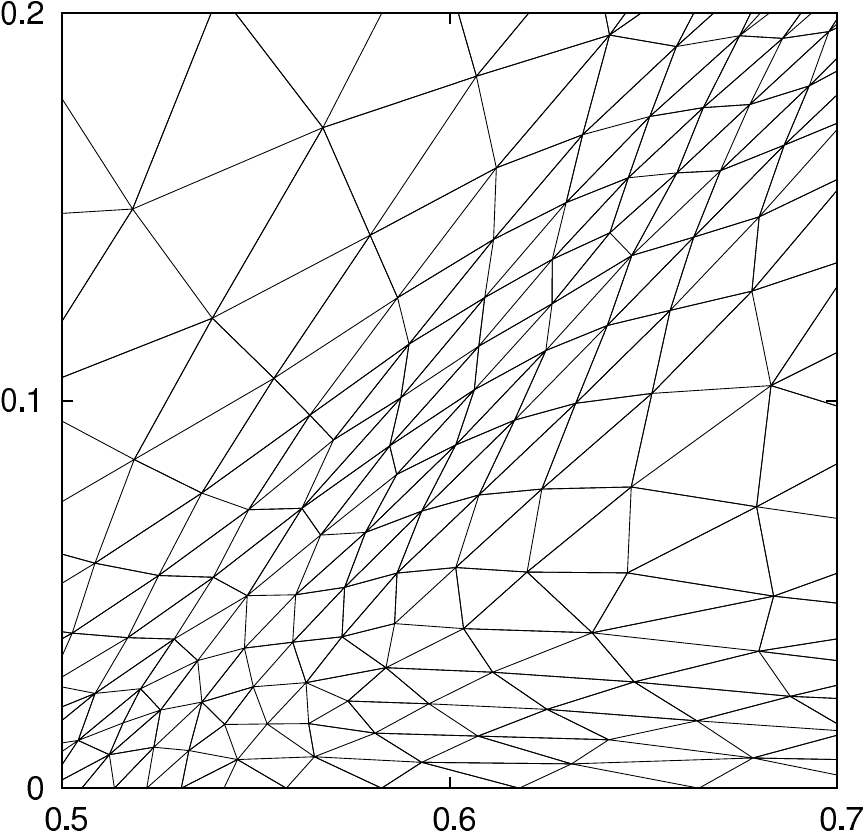}
   }%
   \\[2ex]
   \subcaptionbox{Full error estimator: 
      763 vertices and \num{1413} triangles, 
      $\|e\|_{L^2} = \num{1.6e-3}$, maximum aspect ratio $46.9$.\label{fig:tanh_eh_meshes_global}}[1.0\linewidth]{
      \includegraphics[height=0.28\textheight,clip]{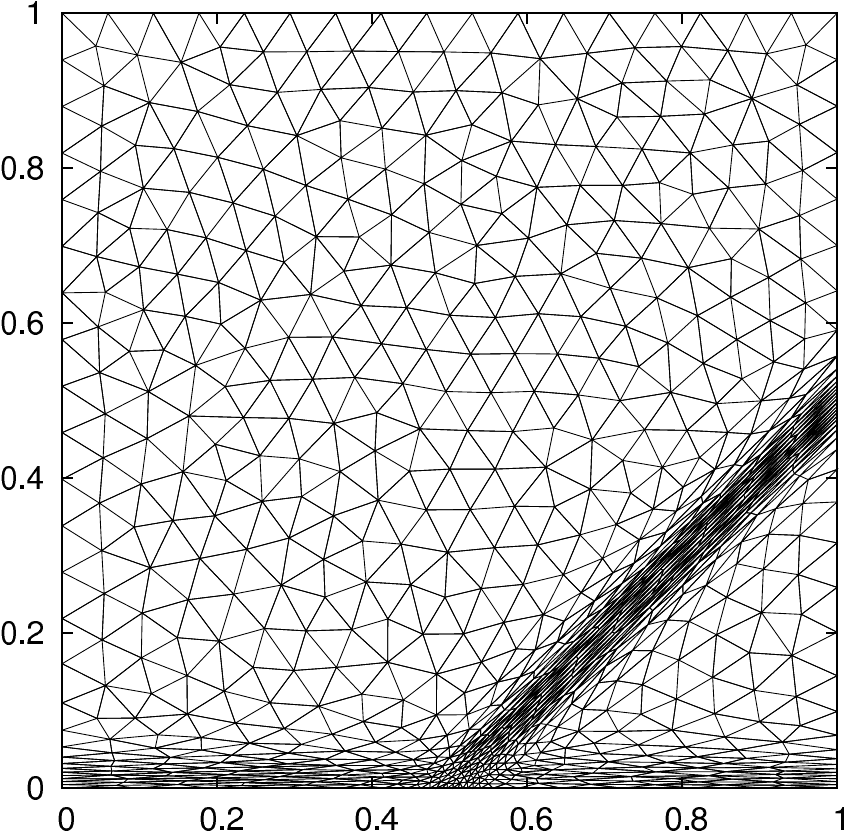}
      \qquad\qquad
      \includegraphics[height=0.28\textheight,clip]{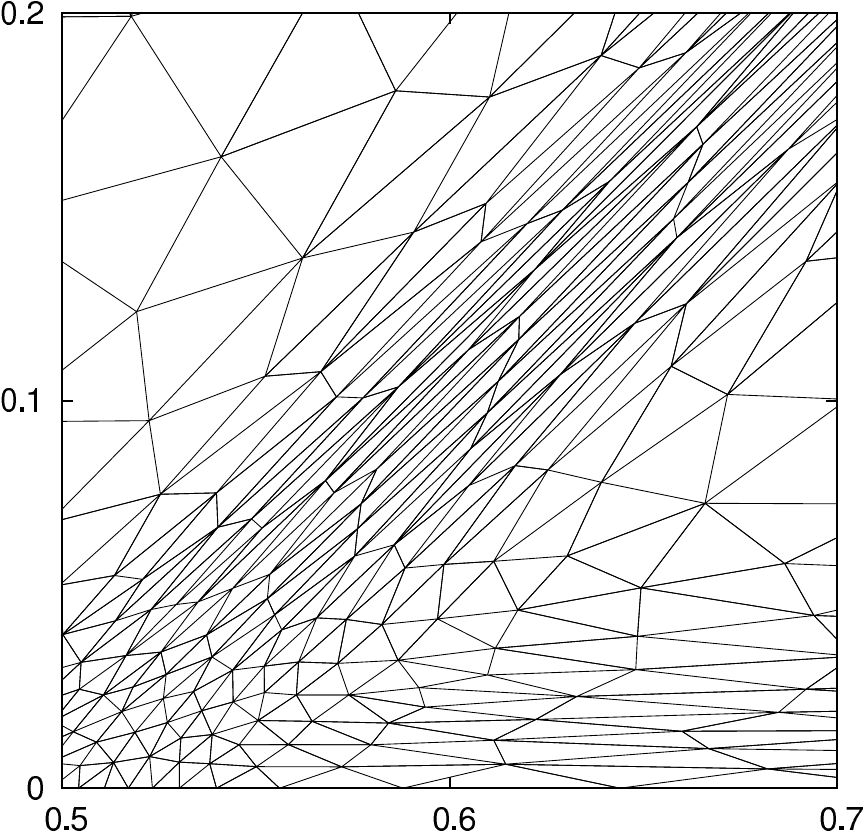}
   }%
   \caption{Example~\ref{ex:tanh}: adaptive meshes obtained by means of the 
      reduced and full a posteriori error estimators (left) 
      and close-up views at $(0.6,0.1)$ (right).
      The desired number of mesh elements is 600.}\label{fig:tanh_eh_meshes}
\end{figure}

Results for the error of the adaptive solution against the number of elements are presented in \cref{fig:tanh_eh_comparison}. 

As expected, the full error estimator works best, leading to a smaller error than those obtained with local error estimators.
The node-based error estimator works better than the edge-based error estimator, mainly because it involves more elements and, in this sense, is more global.

The same observation can also be made from \cref{fig:tanh_eh_meshes}, where adaptive meshes obtained with the error estimators are shown. 
For these mesh example, the desired number of mesh elements $N$ in the normalized metric tensor given by \cref{eq:optimal_metric2+1} has been set to $600$.
All methods produce correct mesh concentrations, although mesh alignment and orientation are different. 
In the mesh controlled by the full error estimator elements near the boundary layer and the shock wave are very thin, have a large aspect ratio\footnote{Aspect ratio is longest edge divided by shortest altitude. An equilateral triangle has an aspect ratio of $\sqrt{3} / 2 \approx 0.9$.} of up to $46.9$, and are properly aligned with the fronts of the shock wave and the boundary layer (\cref{fig:tanh_eh_meshes_global}). 
On the other hand, the elements of meshes controlled by reduced error estimators have rather moderate aspect ratios of $12.8$ and $14.3$ and are less anisotropic (\cref{fig:tanh_eh_meshes_diagonal,fig:tanh_eh_meshes_patch}). 

The accuracy of the corresponding finite element solutions is different, too.
The mesh controlled by the full error estimator leads to a solution error $\|e\|_{L^2} = \num{1.6e-3}$, less then one half of $\|e\|_{L^2} = \num{3.7e-3}$, the error obtained using the node-based error estimator, and about one third of $\|e\|_{L^2} = \num{5.0e-3}$, the error achieved with the edge-based error estimator. 

These results are in good agreement with the comments made in \cref{sec:Computation} that the full error estimator will do a better job than reduced ones for anisotropic mesh adaptation.
Reduced error estimators are able to capture the distribution of the magnitude of the true error and yield a good mesh concentration. 
However, they fail to produce proper mesh alignment, i.e., they does not contain enough information for proper shape and orientation adaptation.

\begin{figure}[t]%
   \subcaptionbox{Effects of the number of Gauß-Seidel iterations used in the 
      solution of the linear system resulting from the error problem 
   $(E_h)$.\label{fig:tanh_gs_steps}}[0.47\linewidth]{
   \includegraphics[width=0.47\textwidth,clip]{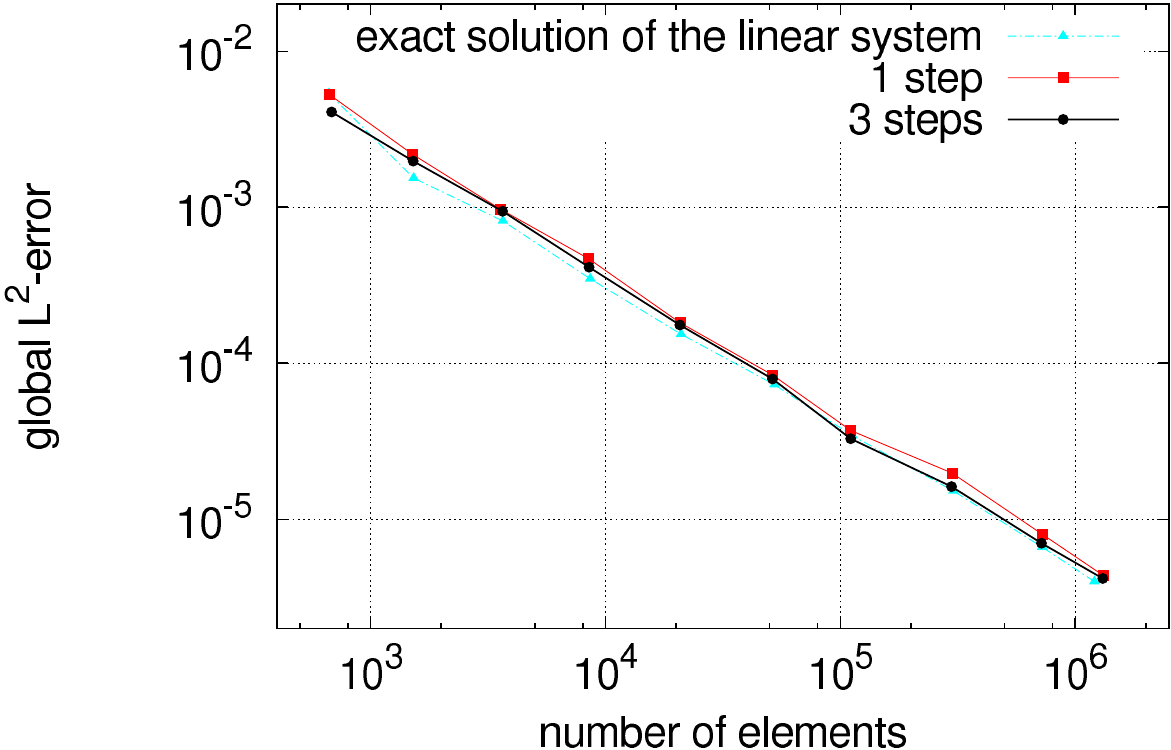}
   }%
   \hfill{}%
   \subcaptionbox{Error comparison for adaptive finite element solutions obtained 
      with global error estimation and Hessian recovery.\label{fig:tanh_wf_qls_hb}}[0.47\linewidth]{
      \includegraphics[width=0.47\textwidth,clip]{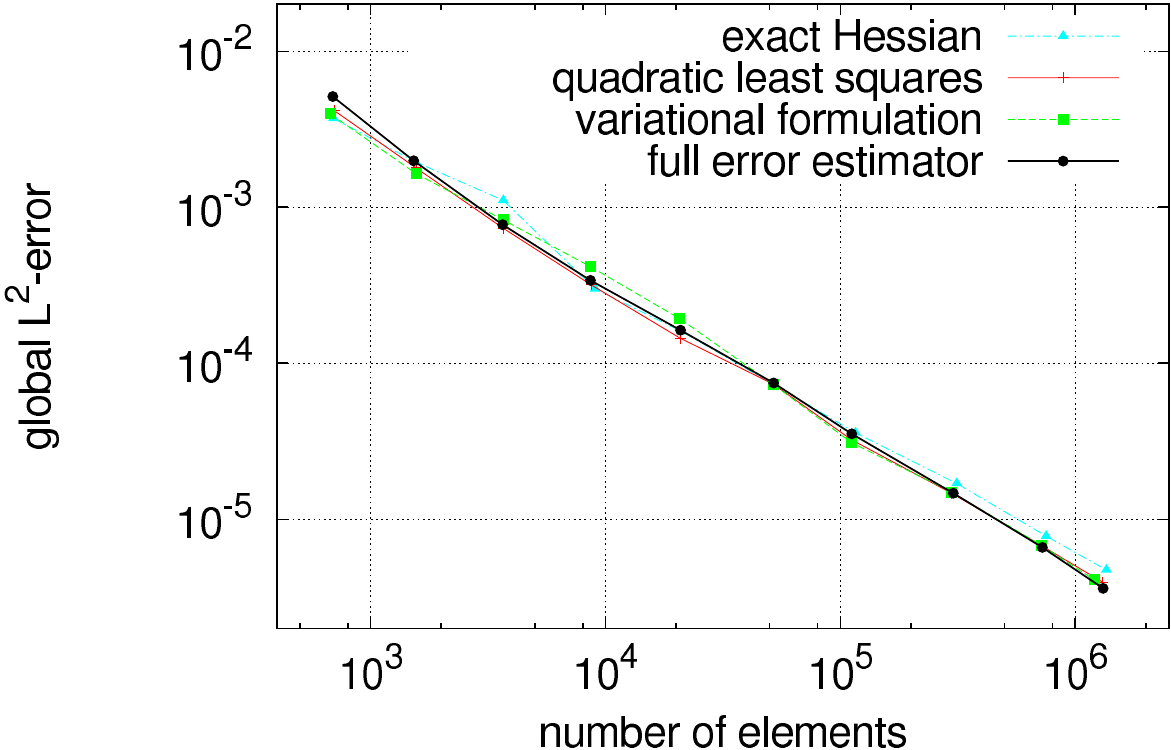}
   }%
   \caption{Example~\ref{ex:tanh}:
      (\protect\subref{fig:tanh_gs_steps}) effects of the number of Gauß-Seidel iterations and 
      (\protect\subref{fig:tanh_wf_qls_hb}) comparison of global error estimation and Hessian recovery.
   }%
   \label{fig:tanh_gs_stepsFig}%
\end{figure}

\textbf{Effects of the number of Gauß--Seidel iterations.}
We now investigate how many iterations are sufficient for obtaining a valuable approximation to the error equation. 
\Cref{fig:tanh_gs_steps} presents results for different iteration numbers to compute the full error estimator. 
As one can see, a few iterations are sufficient for obtaining an approximation good enough for mesh adaptation. 
The convergence lines are very close to each other.
The exact solution of the error problem leads to a smaller error, but the difference is hardly visible.
Three steps of the symmetric Gauß--Seidel method produce an almost optimal mesh for this example. 

\textbf{Comparison to Hessian recovery methods.}
Two Hessian recovery methods are considered for comparison purpose.

\emph{Quadratic least squares fitting.}
This method was recently developed by Zhang and Naga~\cite{ZhaNag05} and proved to be robust and reliable. 
It computes a local quadratic fitting to function values or their approximations at some neighboring points and
obtains a Hessian approximation by differentiating the polynomial twice.

\emph{Variational formulation.} This approach recovers the Hessian, which  does not exists in the classical sense for piecewise linear functions, by means of a variational formulation~\cite{Dolejs98}.
Precisely, let $\phi_i \in V_h$ be the piecewise linear basis function at node $(x_i,y_i)$. Then the nodal approximation to the second-order derivative $u_{xx}$ of a function $u$ at $(x_i,y_i)$ is defined as
\[
\left. (D_{xx}^2u_h)\right |_{(x_i,y_i)}  \int_\Omega \phi_i \,dx\,dy \approx
\int_\Omega D_{xx}^2u_h \,  \phi_i \,dx\,dy = 
   -\int_\Omega \frac{\partial u_h}{\partial x} \,
      \frac{\partial \phi_i }{\partial x} \, dx \, dy.
\]
The same approach is used to approximate $u_{xy}$ and $u_{yy}$.

\Cref{fig:tanh_wf_qls_hb} shows the error against the number of elements for each method.
For comparison purpose, results obtained using the analytical Hessian are also included.
All methods provide almost the same results. Particularly, the method based on the global estimator with three Gauß--Seidel iterations is comparable to the recovery-based methods.

It is worth noting that although the quadratic least squares fitting is generally more accurate and robust than the variational method, both produce basically the same adaptive mesh. 
This seems to confirm the conjecture that highly accurate Hessian recovery is not necessary
for good mesh adaptation.

\subsection{Further examples}\label{ex:problem4}
We consider two boundary value problems in the form \cref{eq:bvp} with now
the right-hand side $f$ and the Dirichlet boundary condition being chosen such that the exact solution
is given by the following functions:
\begin{align*}
   u_1(x,y) &= \frac{1}{1 + e^{\frac{x+y-1.25}{0.05}}}, \\[0.3em]
   u_2(x,y) &=e^{-25x} + e^{-25y}.
\end{align*}
The first function represents a shock wave along the line $y = 1.25-x$ while the second models
a boundary layer near the coordinate axes.

\begin{figure}[p]\centering{}%
   \subcaptionbox{Error comparison for adaptive solutions.}{
      \includegraphics[height=0.185\textheight,clip]{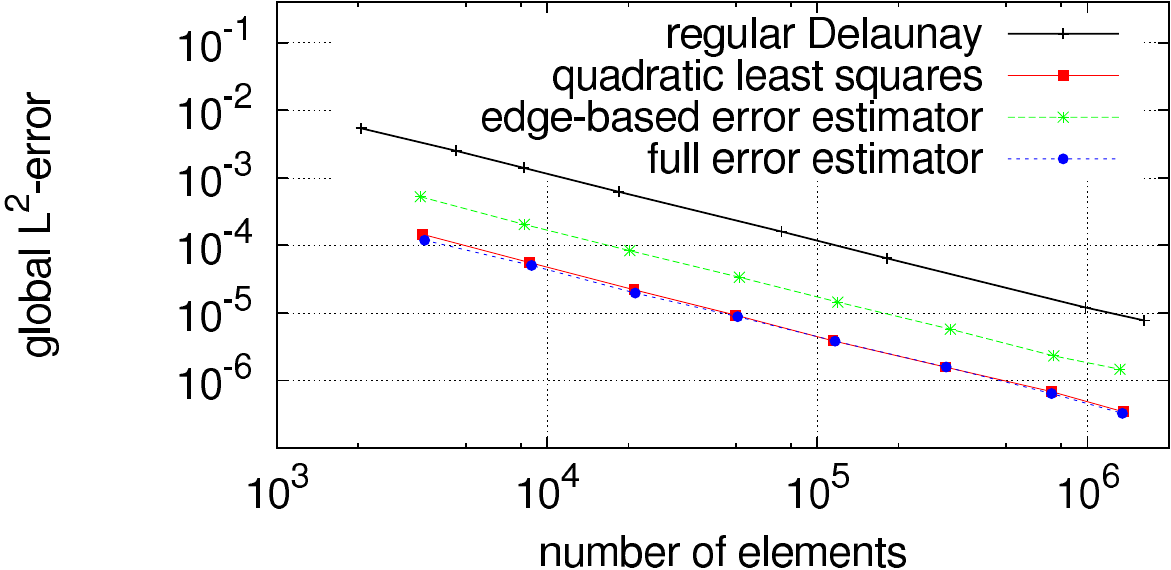}
   }%
   \hfill{}%
   \subcaptionbox{Error estimator effectivity index and $\beta$.\label{fig:p4_effIndex}}{
      \includegraphics[height=0.185\textheight,clip]{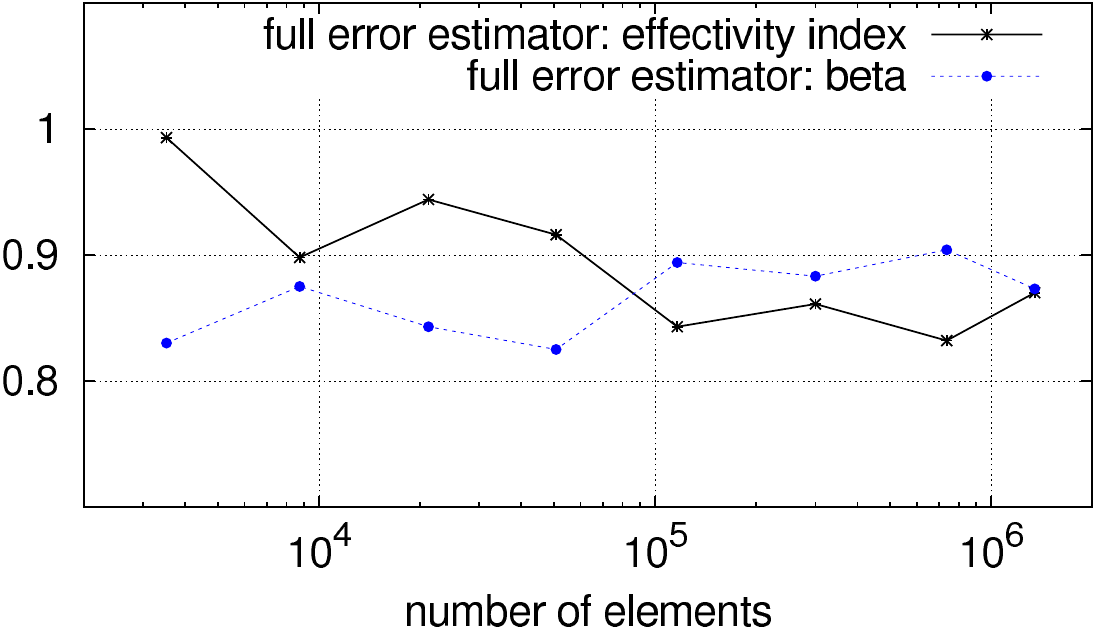}
   }%
   \\[2ex]
   \subcaptionbox{Edge-based error estimator: 
      $684$ vertices and \num{1282} triangles,
      $\|e\|_{L^2} = \num{1.4e-3}$, maximum aspect ratio $3.8$.}[0.9\linewidth]{
      \includegraphics[height=0.20\textheight,clip]{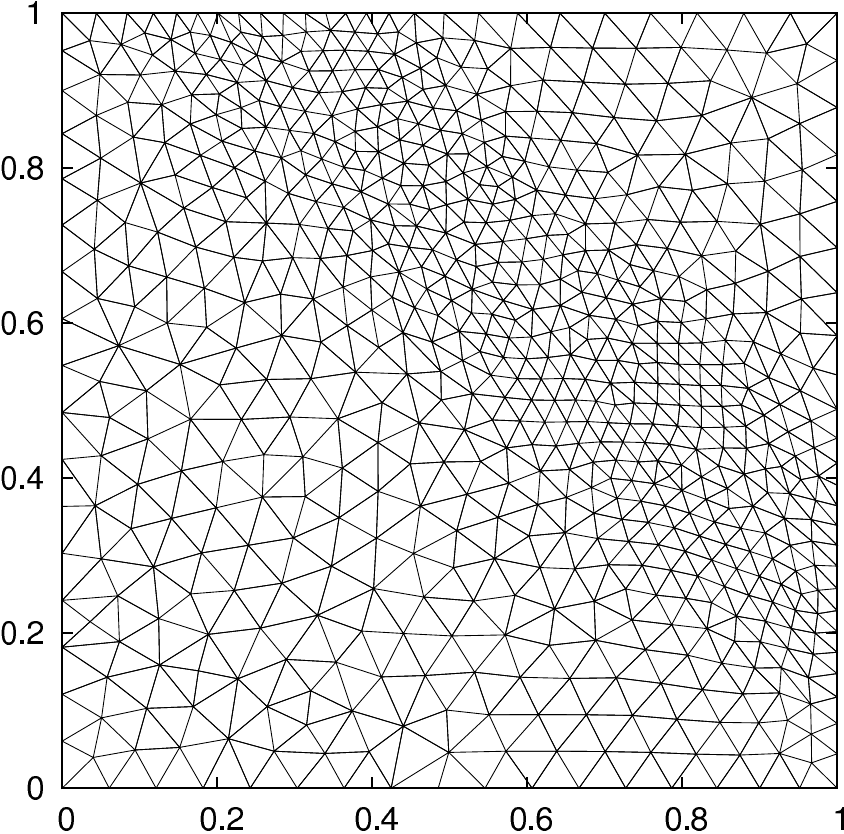}
      \qquad
      \includegraphics[height=0.20\textheight,clip]{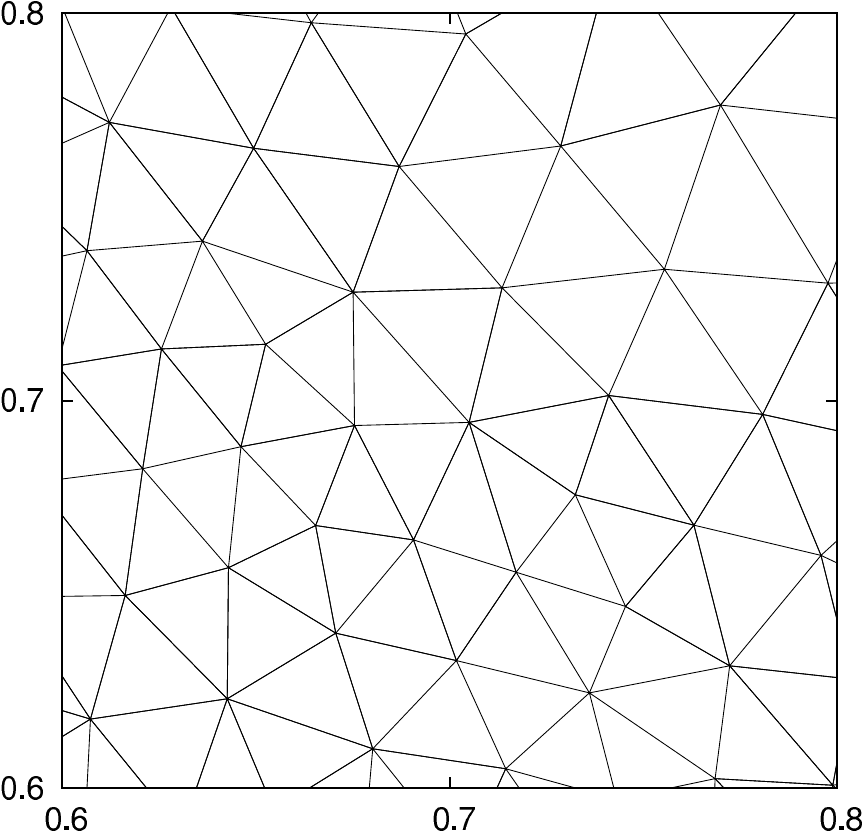}
   }%
   \\[2ex]
   \subcaptionbox{Quadratic least squares Hessian recovery: 
      $693$ vertices and \num{1277} triangles, $\|e\|_{L^2} = \num{3.5e-4}$,
      maximum aspect ratio $13.5$.}[0.95\linewidth]{
      \includegraphics[height=0.20\textheight,clip]{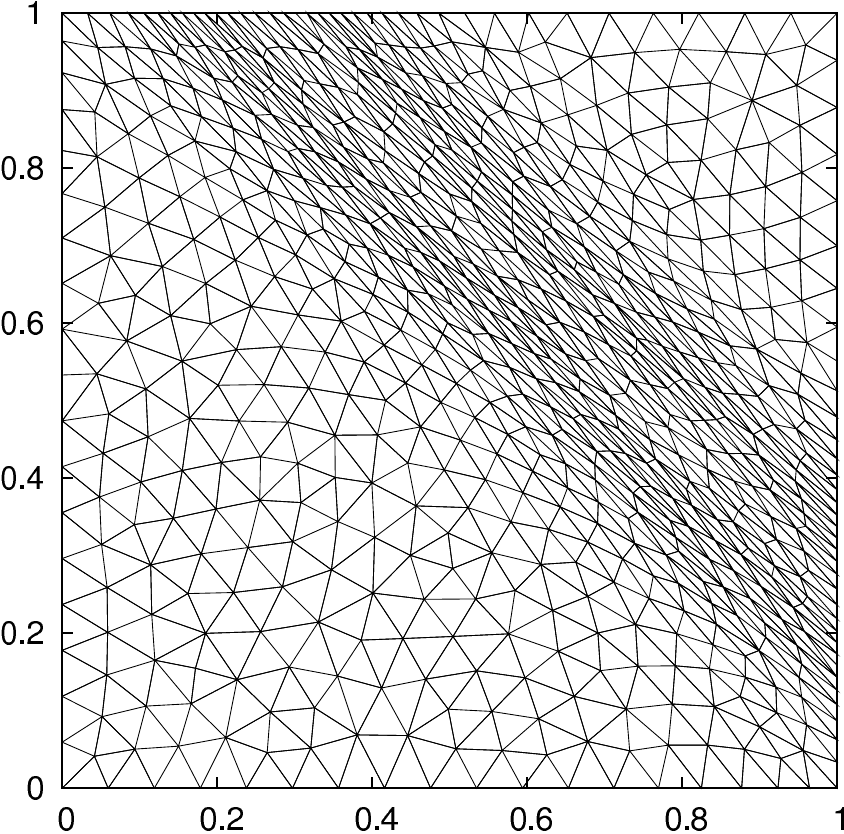}
      \qquad
      \includegraphics[height=0.20\textheight,clip]{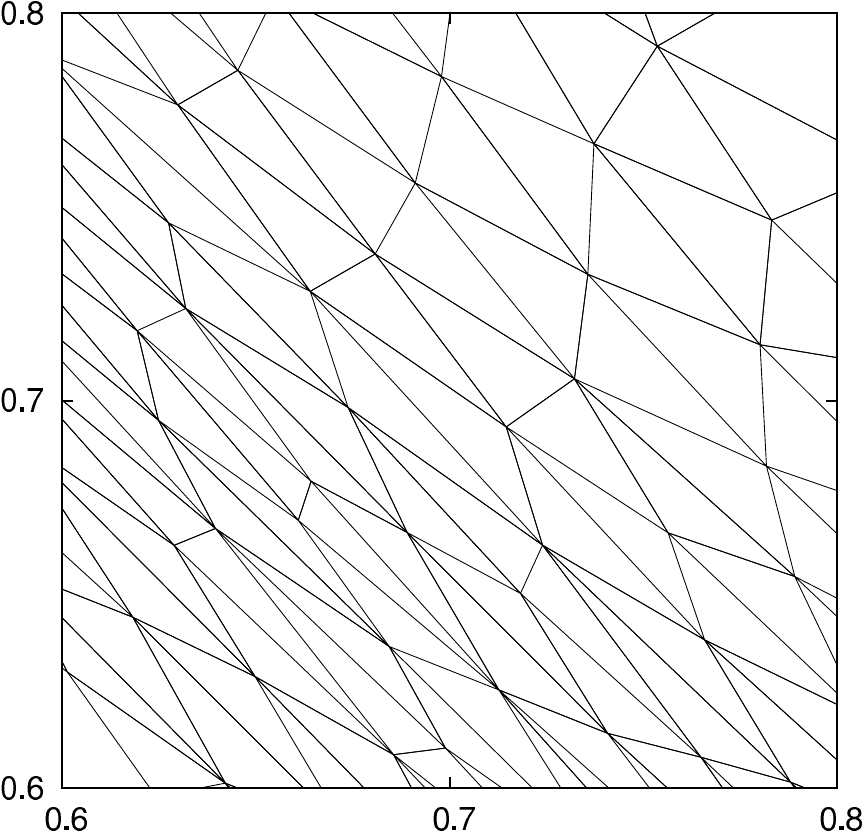}
   }%
   \\[2ex]
   \subcaptionbox{Full error estimator: 
      $714$ vertices and \num{1318} triangles, $\|e\|_{L^2} = \num{3.4e-4}$,
      maximum aspect ratio $15.0$.}[0.9\linewidth]{
      \includegraphics[height=0.20\textheight,clip]{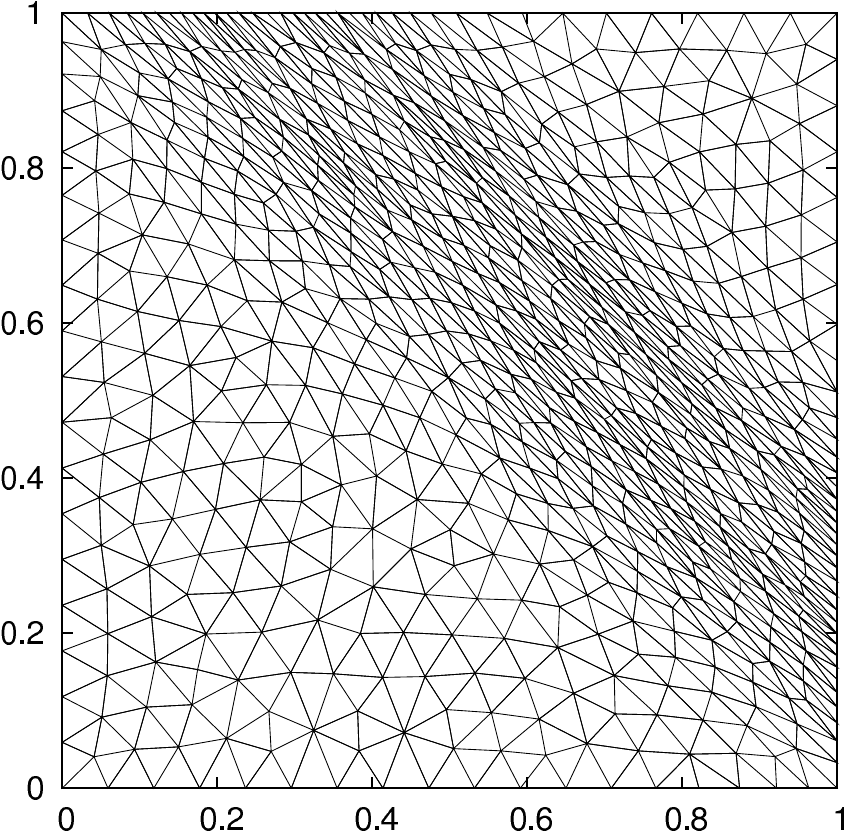}
      \qquad
      \includegraphics[height=0.20\textheight,clip]{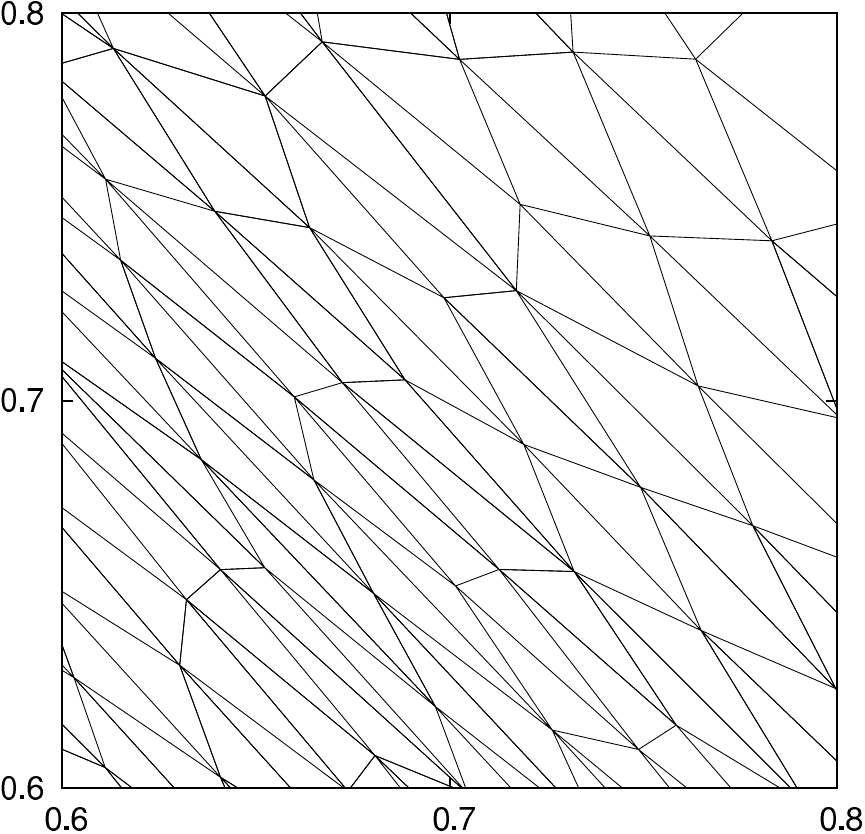}
   }%
   \caption{BVP \cref{eq:bvp} with the exact solution $u(x,y) = 1/(1 + e^{\frac{x+y-1.25}{0.05}})$: 
      adaptive meshes (left) and close-up views at $(0.7,0.7)$ (right).
      The desired number of mesh elements is $600$.}\label{fig:problem4_qls_hb_uniform}
\end{figure}

\begin{figure}[p] \centering
   \subcaptionbox{Error comparison for adaptive solutions.}{
      \includegraphics[height=0.185\textheight,clip]{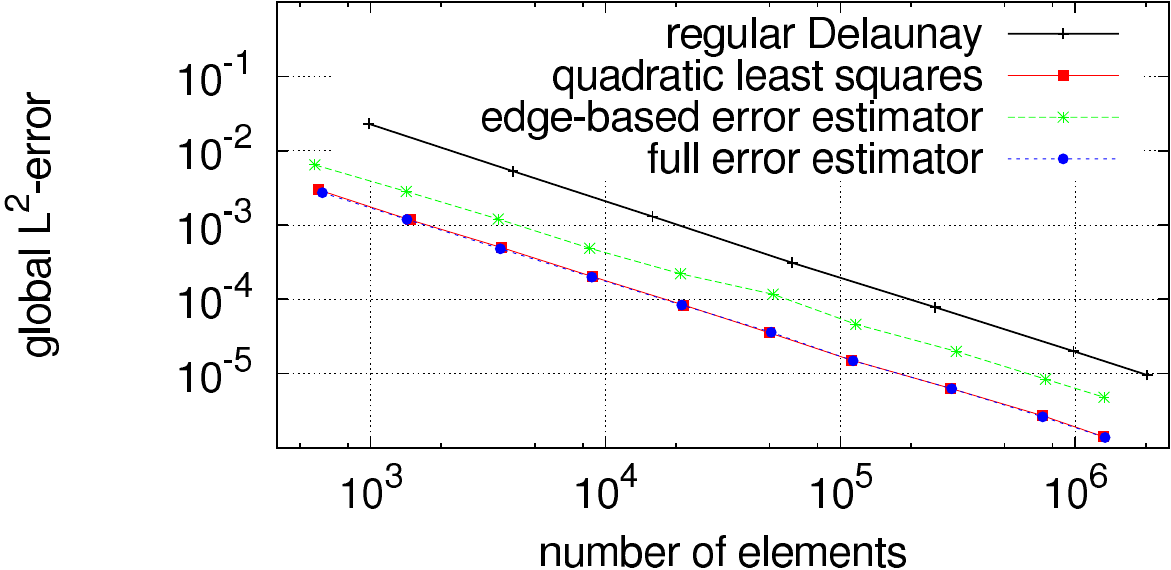}
   }\hfill{}%
   \subcaptionbox{Error estimator effectivity index and $\beta$.\label{fig:bl_effIndex}}{
      \includegraphics[height=0.185\textheight,clip]{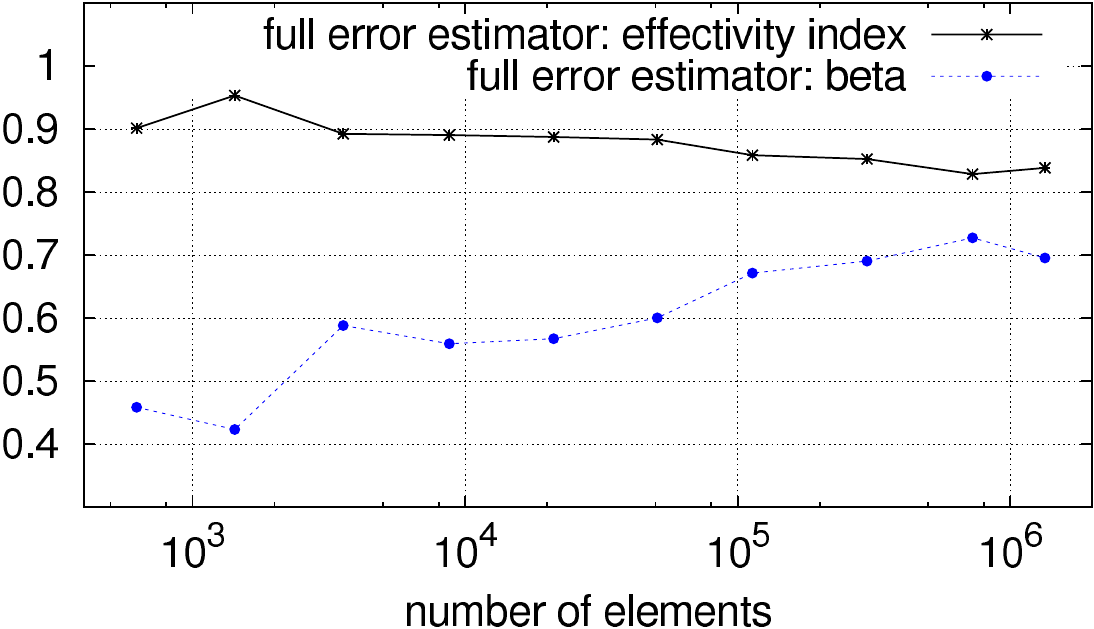}
   }\\[2ex]
   \subcaptionbox{Edge-based error estimator: 
      $698$ vertices and \num{1298} triangles, $\|e\|_{L^2}=\num{3.0e-3}$,
      maximum aspect ratio $7.0$.}[0.9\linewidth]{
      \includegraphics[height=0.20\textheight,clip]{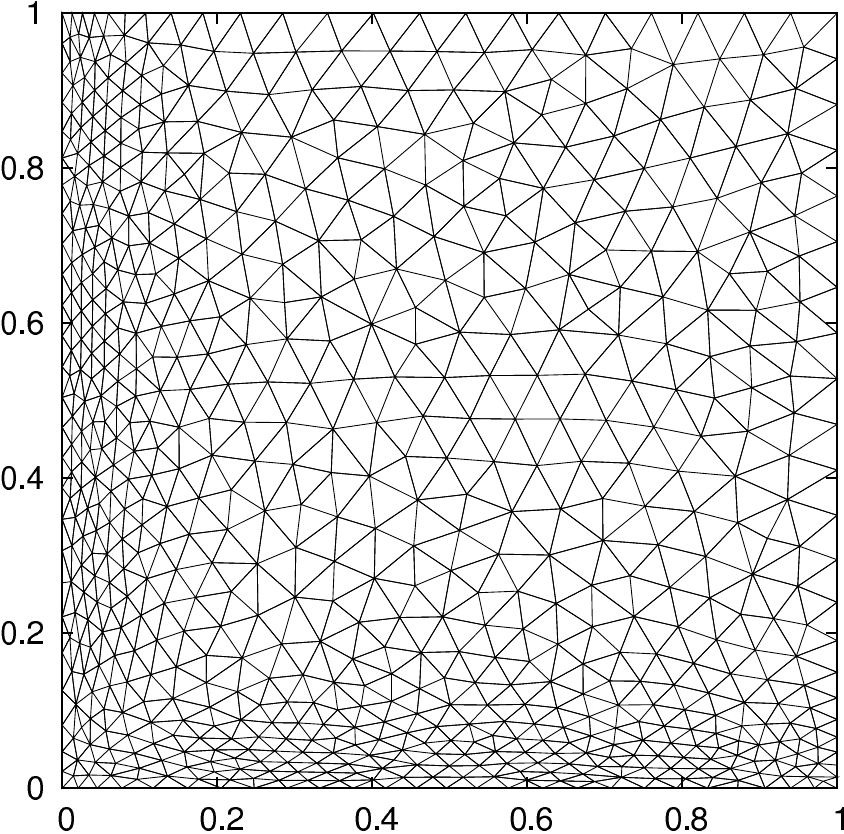}
      \qquad
      \includegraphics[height=0.20\textheight,clip]{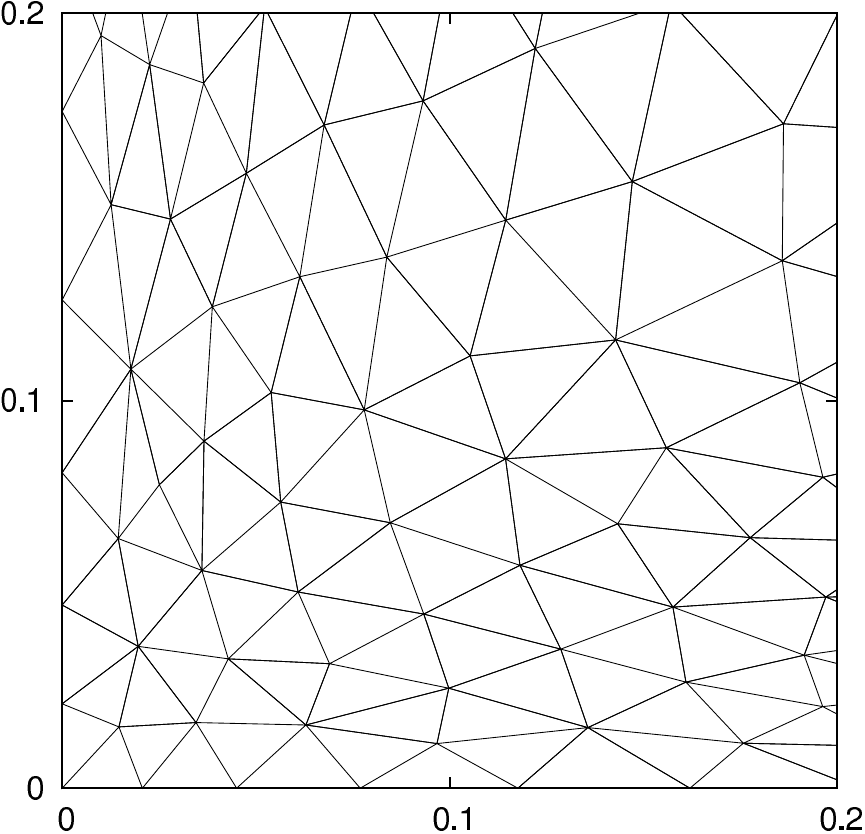}
   }%
   \\[2ex]
   \subcaptionbox{Quadratic least squares Hessian recovery:
      $710$ vertices and \num{1327} triangles, $\|e\|_{L^2}=\num{1.3e-3}$,
      maximum aspect ratio $20.0$.}[0.95\linewidth]{
      \includegraphics[height=0.20\textheight,clip]{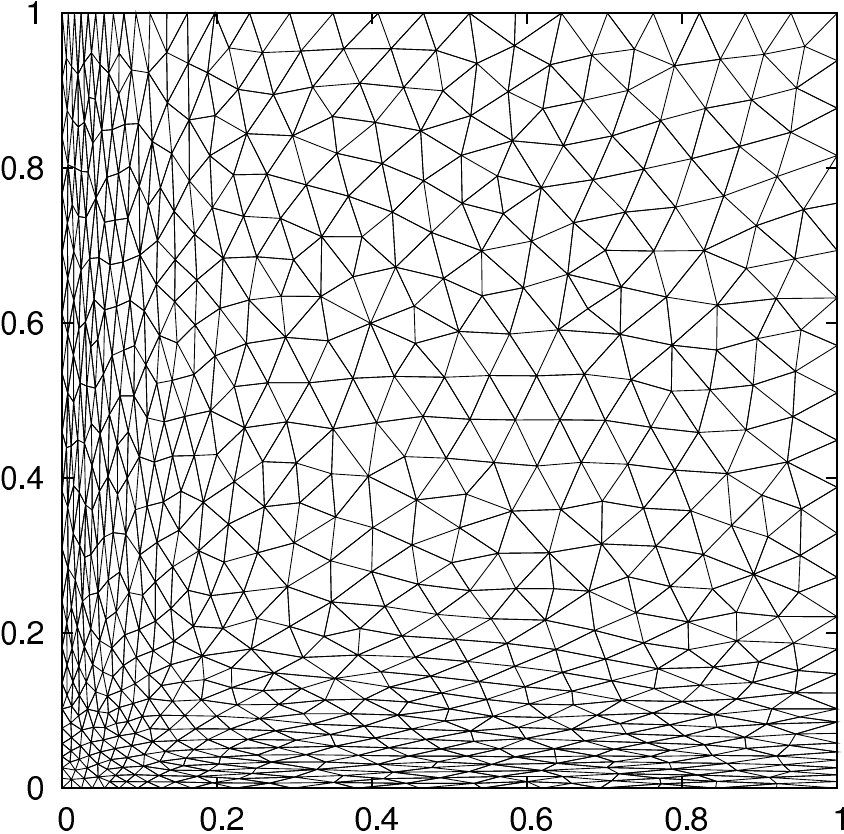}
      \qquad
      \includegraphics[height=0.20\textheight,clip]{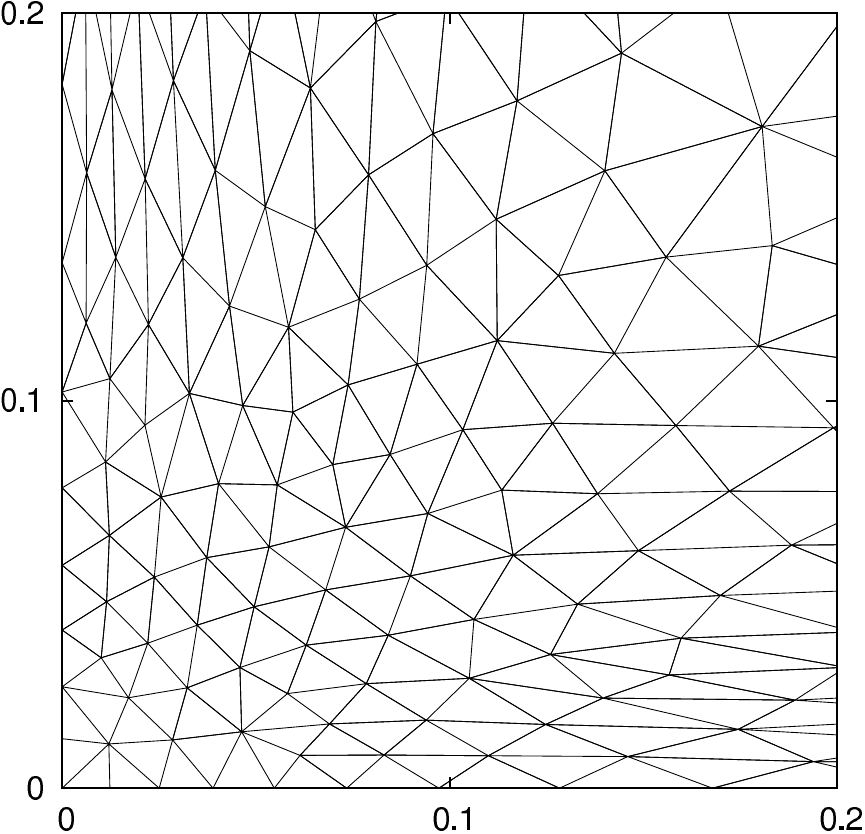}
   }%
   \\[2ex]
   \subcaptionbox{Full error estimator:
      $712$ vertices and \num{1332} triangles, $\|e\|_{L^2}=\num{1.3e-3}$,
      maximum aspect ratio $17.8$.}[0.9\linewidth]{
      \includegraphics[height=0.20\textheight,clip]{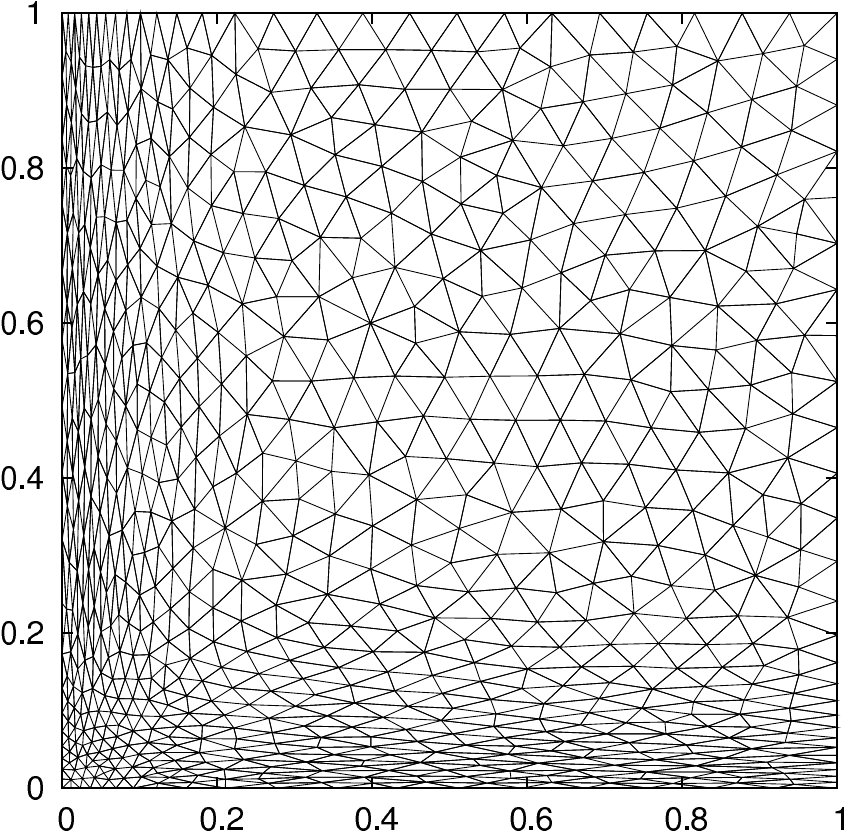}
      \qquad
      \includegraphics[height=0.20\textheight,clip]{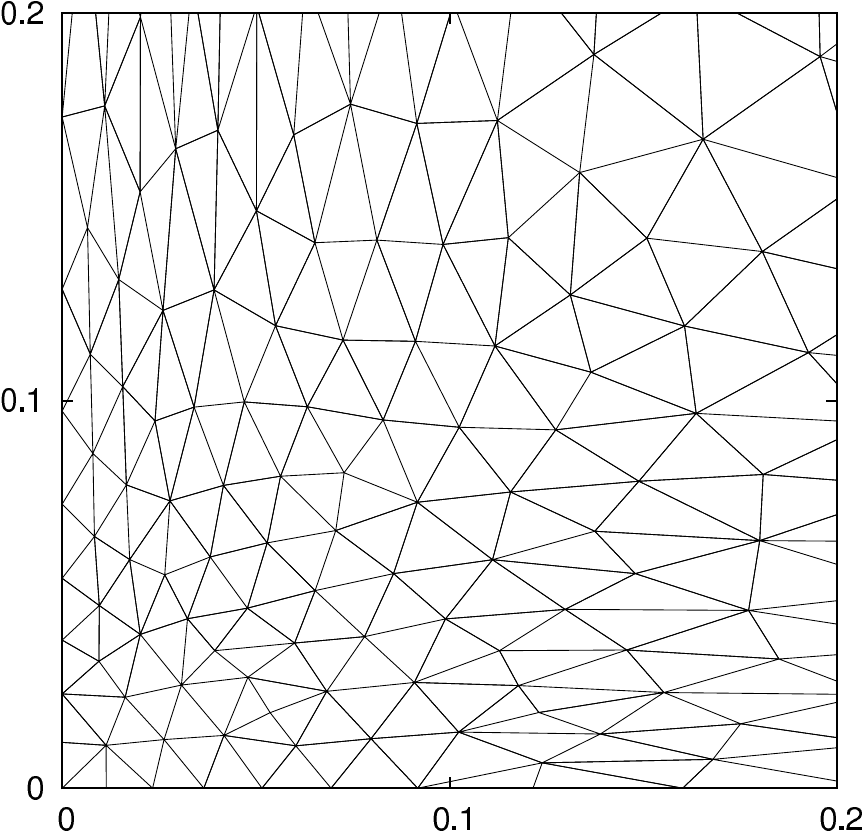}
   }%
   \caption{BVP~\eqref{eq:bvp} with the exact solution 
      $u(x,y) = e^{-25x} + e^{-25y}$: 
      adaptive meshes (left) and close-up views at $(0.1,0.1)$ (right).
      The desired number of mesh elements is 600.}\label{fig:bl_qls_hb_uniform}
\end{figure}

We compare the error for finite element solutions obtained with the global error estimator and the quadratic least squares Hessian recovery. 
Results for the quasi-uniform (regular Delaunay) mesh and the edge-based error estimator are also given.
\Cref{fig:problem4_qls_hb_uniform,fig:bl_qls_hb_uniform} show the results.

As in \cref{ex:tanh}, we can see that mesh adaptation significantly reduces the finite element error compared to a quasi-uniform mesh having the same number of elements. The mesh based on the edge-based error estimator provides a good mesh concentration and is clearly better than a quasi-uniform one, but it is almost isotropic and inferior to a mesh obtained with the use of the full error estimator.
Again, one can observe that the elements of the meshes obtained by means of the full error estimator and the quadratic least squares fitting are properly aligned with the shock wave and the boundary layers. 
Thus, the new method produces results comparable to those obtained with recovery-based methods.

\subsection{Discontinuous gradients}\label{sec:jG}
\begin{figure}[hp] \centering
   \subcaptionbox{
      Side views of finite element solutions: 
      without the interface being present in the mesh, the solution is 
      non-planar and not exact (left, different lines represent the solution
      for different values of $y$); 
      with the interface being present in the mesh, the solution is exact 
   (right).\label{fig:jg_surface_plots}}[0.95\linewidth]{
      \includegraphics[width=0.35\textwidth,clip]{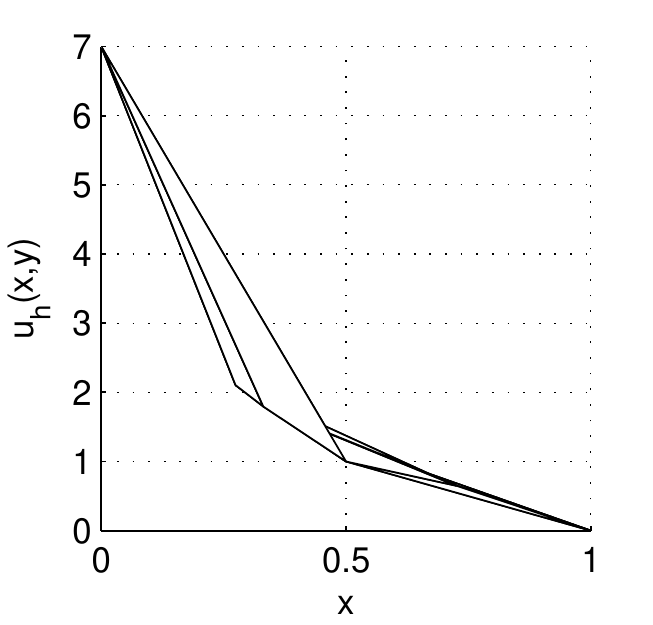}
      \qquad\qquad
      \includegraphics[width=0.35\textwidth,clip]{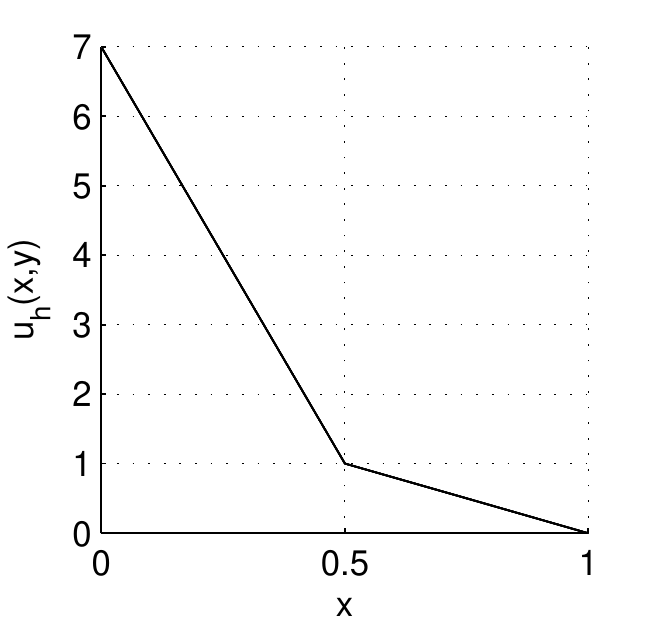}
   }%
   \\[2ex] 
   \subcaptionbox{Adaptation without predefined interface edges: 
      quadratic least squares Hessian recovery, $385$ vertices, 
      $\|e\|_{L^2}=\num{4.2e-3}$, 
      maximum aspect ratio $8.1$ (left); 
      full error estimator, $381$ vertices, 
      $\|e\|_{L^2}=\num{4.1e-3}$,
      maximum aspect ratio $7.1$ (right).\label{fig:jg_without_interface}}[0.95\linewidth]{
      \includegraphics[width=0.35\textwidth,clip]{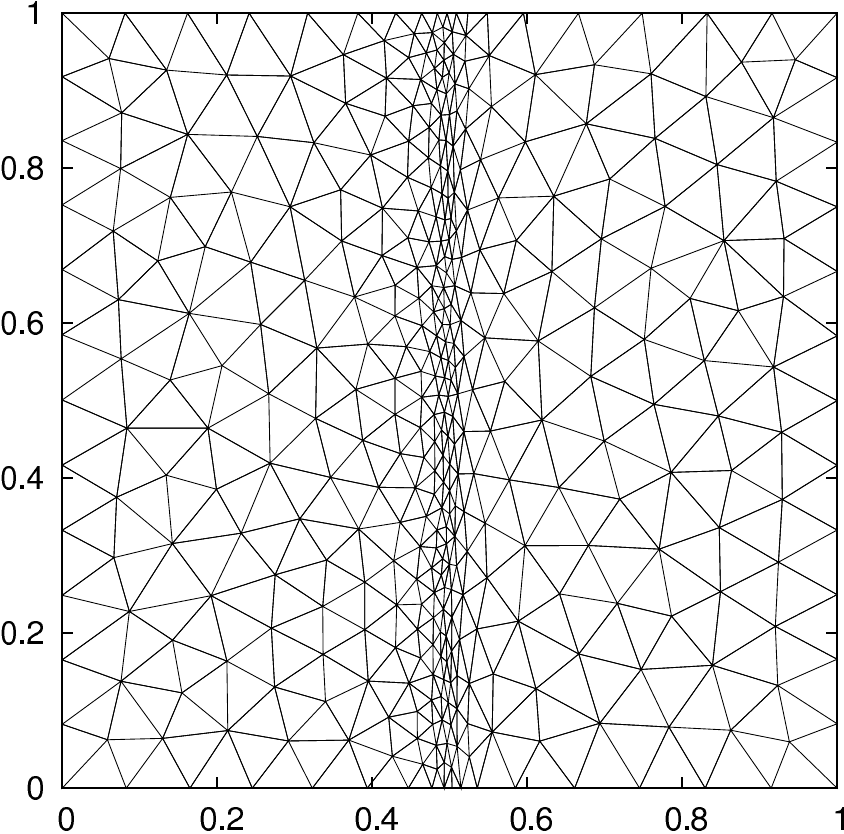}
      \qquad\qquad
      \includegraphics[width=0.35\textwidth,clip]{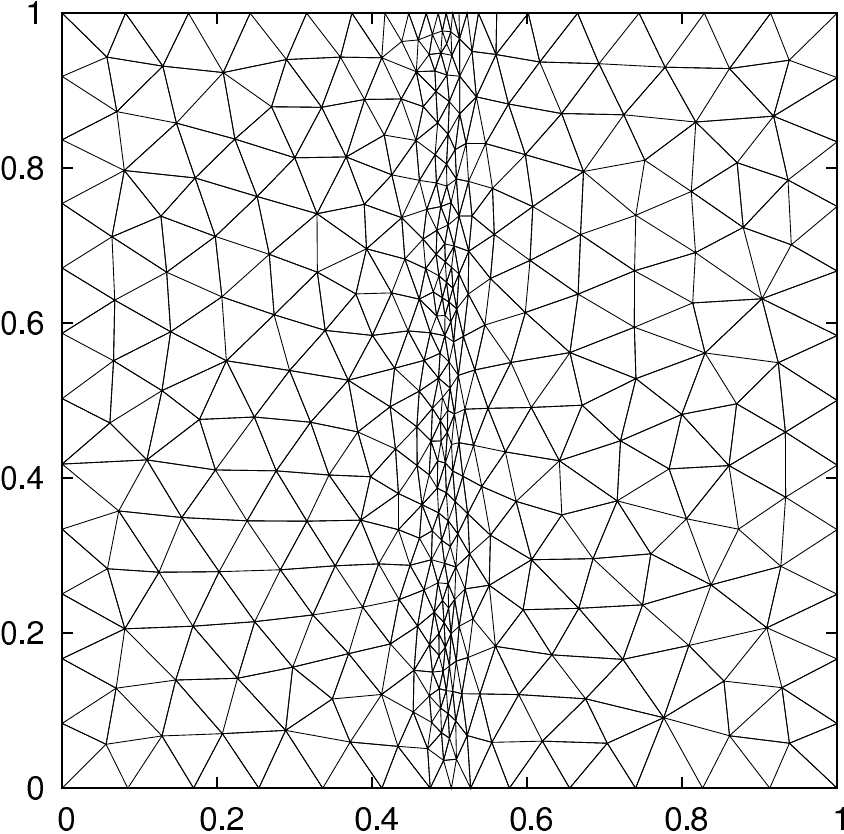}
   }%
   \\[2ex]
   \subcaptionbox{Adaptation with predefined interface edges:
      quadratic least square Hessian recovery: $77$ vertices, 
      $\|e\|_{L^2} = \num{3.1e-16}$, 
      maximum aspect ratio $60.7$ (left);      
      full error estimator, $63$ vertices, 
      $\|e\|_{L^2} = \num{3.2e-16}$,
      maximum aspect ratio $1.8$ (right).\label{fig:jg_with_interface}}[0.95\linewidth]{
      \includegraphics[width=0.35\textwidth,clip]{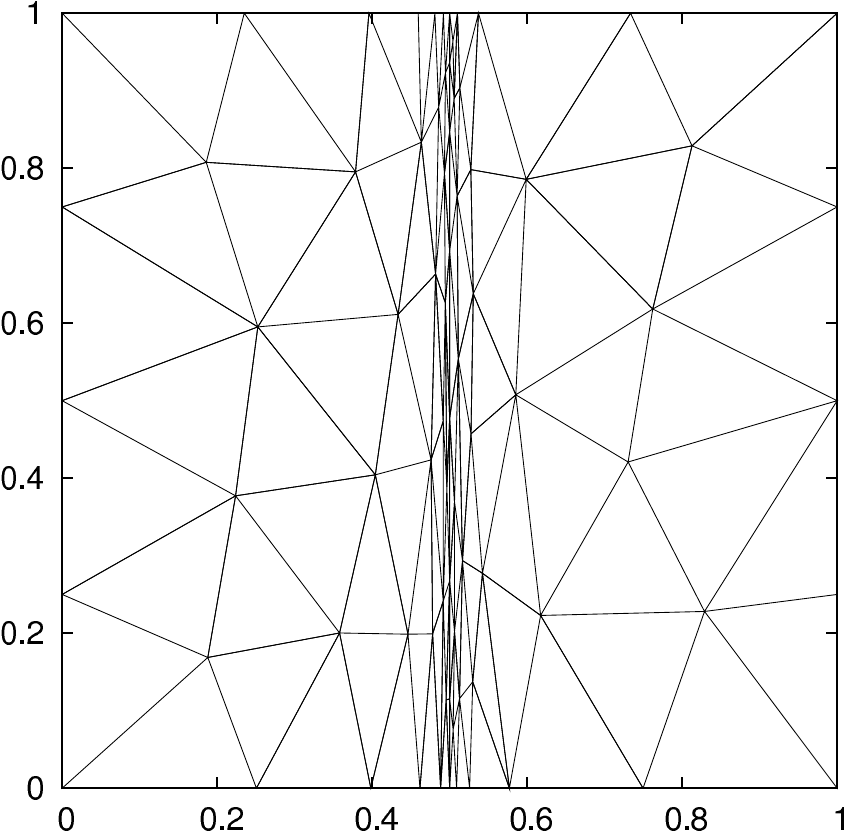} 
      \qquad\qquad
      \includegraphics[width=0.35\textwidth,clip]{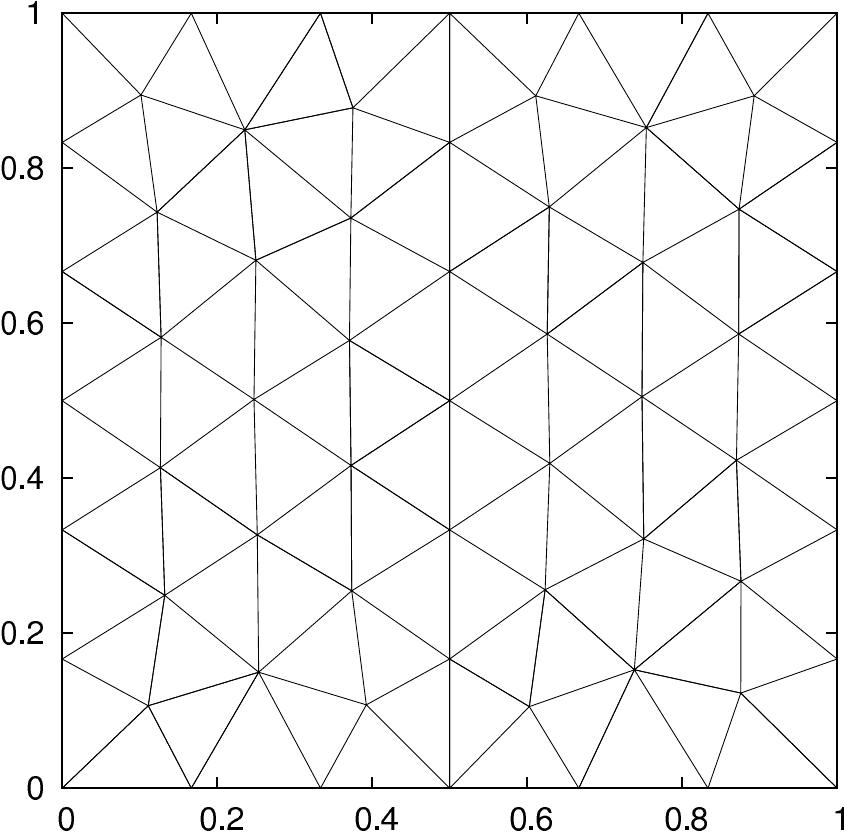}
   }%
   \caption{Example~\ref{sec:jG}: gradient jump along the 
      line $x=0.5$. Adaptive meshes and finite element solutions with and 
      without the predefined interface edges.}\label{fig:jg}
 \end{figure}
Next, we consider problems whose solution has a discontinuous gradient along a certain interface in the domain. This situation arises in elliptic problems with discontinuous coefficients in the diffusion term such as heat conduction
problems with jumps in material coefficients. Difficulties when using gradient recovery methods for such problems were
already pointed out in~\cite{Ovall06a}, and this is true for the Hessian recovery as well: if the numerical approximation is accurate enough, we should expect a discontinuity in its gradient and its Hessian. Since most Hessian recovery methods employ some sort of averaging over a certain region, they can be very inaccurate near discontinuities.
This issue can readily be observed in the following simple example.

Let $\Omega = (0,1) \times (0,1)$. Consider the boundary value problem 
\[
   \begin{cases}
   - a \Delta  u = 0  &\text{in } \Omega, \\
               u = g  &\text{on } \partial \Omega,
   \end{cases} 
\]
where
\[ a = \begin{cases}
    1,      & x <  0.5, \\ 
    \alpha, & x \geq 0.5
   \end{cases} \]
and the Dirichlet boundary condition is chosen such that the exact solution is given by 
\begin{align*}
   u(x,y) = 
   \begin{cases}
      -2 \alpha x + \alpha + 1,  &x <  0.5, \\
      - 2 x + 2,                 &x \geq 0.5.
   \end{cases}
   \label{ex:jumpingGradient}
\end{align*}
The solution has a gradient jump of magnitude $\alpha$ across the line $x=0.5$, but is continuous on $\Omega$ and
linear in each of the subdomains. We take $\alpha = 6$ in our computation.

We first consider the situation where the mesh does not contain the information of the interface. In this situation
at least part of the interface does not consist of edges. 
In order to match the sharp bend in the solution along the interface, the adaptive mesh should exhibit a strong concentration of elements around $x=0.5$ oriented along the interface. In this test, the quadratic least squares and
the global hierarchical basis error estimator both succeed in providing an appropriate mesh adaptation and, again, deliver comparable results (\cref{fig:jg_without_interface}).

The situation is different if the interface is present in the mesh. In this case, the analytical solution $u$ belongs to the corresponding finite element space and, consequently, the numerical approximation computed by means of the linear finite element method is exact (\cref{fig:jg_surface_plots}, right). Hence, no adaptation is required and the proper mesh should be a uniform mesh. Now, consider the mesh adaptation using the quadratic least squares Hessian recovery. Because of the sharp bend in the solution, the recovered Hessian should be very large, $\approx \mathcal{O}(1/h)$, near $x=0.5$, but zero elsewhere, because the solution is linear in each of the subdomains. This should lead to an excessive over-adaptation near the interface. On the other hand, we expect no adaptation for the hierarchical basis error estimator in this case
because the numerical solution is exact and, consequently, the error estimator is zero everywhere in $\Omega$.
A quasi-uniform mesh should result.
\Cref{fig:jg_with_interface} presents mesh examples. We see that the adaptation by means of the Hessian recovery (left) leads to a strong element concentration along the interface line, as predicted, whereas the mesh based on the hierarchical error estimator (right) is almost uniform.

We also expect a similar behaviour of these methods for general problems exhibiting gradient jumps or similar discontinuities along internal interfaces. Thus, for such problems, it can be of advantage to use the a posteriori error estimator for effective mesh adaptation because of the more efficient employment of given degrees of freedom.
\section{Heat conduction in a thermal battery}\label{sec:Real-Life}
In this section, we consider heat conduction in a thermal battery with large orthotropic jumps in the material coefficients.
The mathematical model considered here is taken from~\cite{Ovall06a,ParDem06} and described by
\begin{equation} 
   \begin{cases} 
      \nabla \cdot ( D^k \nabla u) = f^k       &\text{in } \Omega, \\
      D^k\nabla u \cdot n = g^i - \alpha^i u   &\text{on } \partial\Omega,
   \end{cases} 
   \label{ex:battery}
\end{equation}
where $\Omega = (0,8.4) \times (0,24)$ and
\[ D^k = \begin{bmatrix}
      D_x^k &  0 \\ 
      0     &  D_y^k 
   \end{bmatrix}. \]
The data for each material $k$ and for each of the four sides $i$ of the boundary starting with the left-hand side boundary and ordering them clockwise are given in \cref{tab:battery_data}.
\begin{table}[t]\centering{}%
   \caption{Heat conduction in a thermal battery: material
      coefficients and boundary conditions.}\label{tab:battery_data}%
   \begin{subtable}[t]{.4\linewidth}\centering{}%
      \caption{Material coefficients.}%
      \begin{tabular}{c|ccc}%
         \toprule
         Region $k$ &  $D_x^k$  &  $D_y^k$  &  $f^k$ \\
         \midrule
         1   &  25       &  25       &  0     \\ 
         2   &  7        &  0.8      &  1     \\
         3   &  5        &  0.0001   &  1     \\
         4   &  0.2      &  0.2      &  0     \\
         5   &  0.05     &  0.05     &  0     \\
         \bottomrule
      \end{tabular}%
   \end{subtable}%
   \begin{subtable}[t]{.4\linewidth}\centering{}%
      \caption{Boundary conditions.}%
      \begin{tabular}{c|cc}
         \toprule
         Boundary $i$ &  $\alpha^i$  &  $g^i$ \\
         \midrule
                  1   &  0           &  0    \\ 
                  2   &  1           &  3    \\
                  3   &  2           &  2    \\
                  4   &  3           &  0    \\
         \bottomrule
      \end{tabular}%
   \end{subtable}%
\end{table}%
\begin{figure}[t]\centering{}%
   \subcaptionbox{\label{fig:battery_geometry}}{
      \includegraphics[width=0.2\textwidth, height=0.56\textwidth,clip]{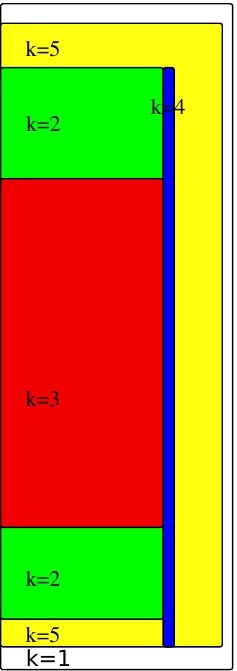}
   }%
   \qquad%
   \subcaptionbox{\label{fig:battery_solution_}}{
      \includegraphics[width=0.2\textwidth, height=0.56\textwidth,clip]{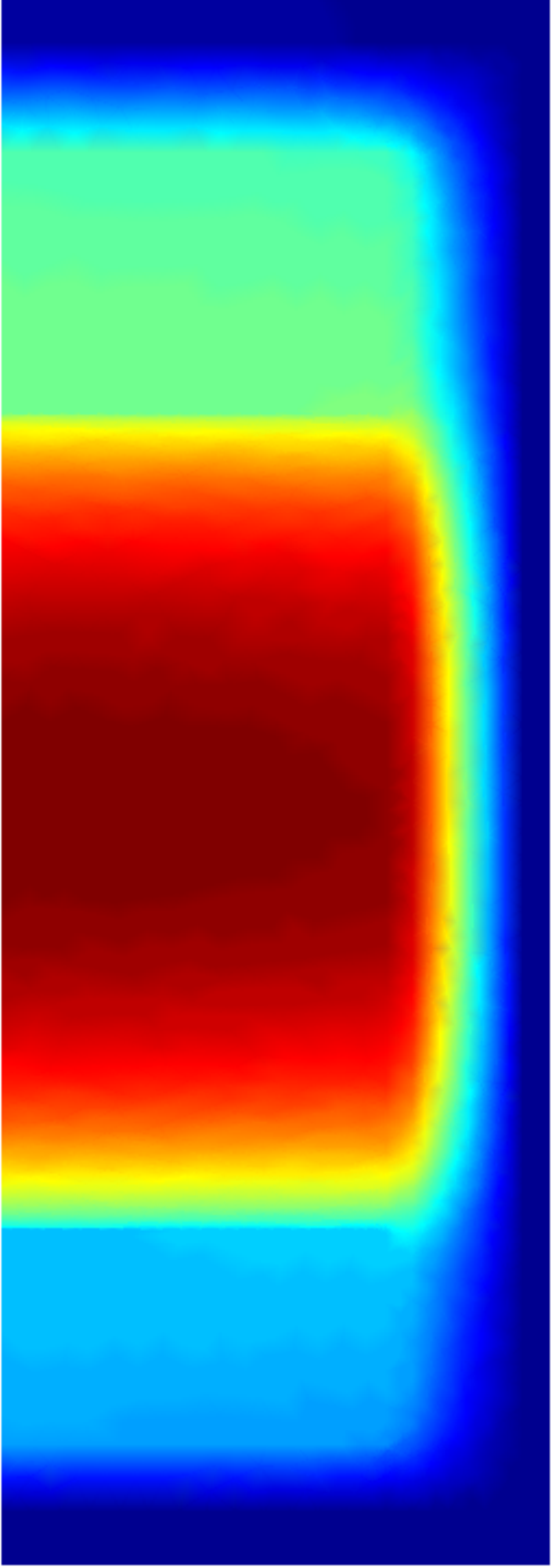}
   }%
   \qquad%
   \subcaptionbox{\label{fig:battery_surface}}{
      \includegraphics[width=0.2\textwidth, height=0.56\textwidth,clip]{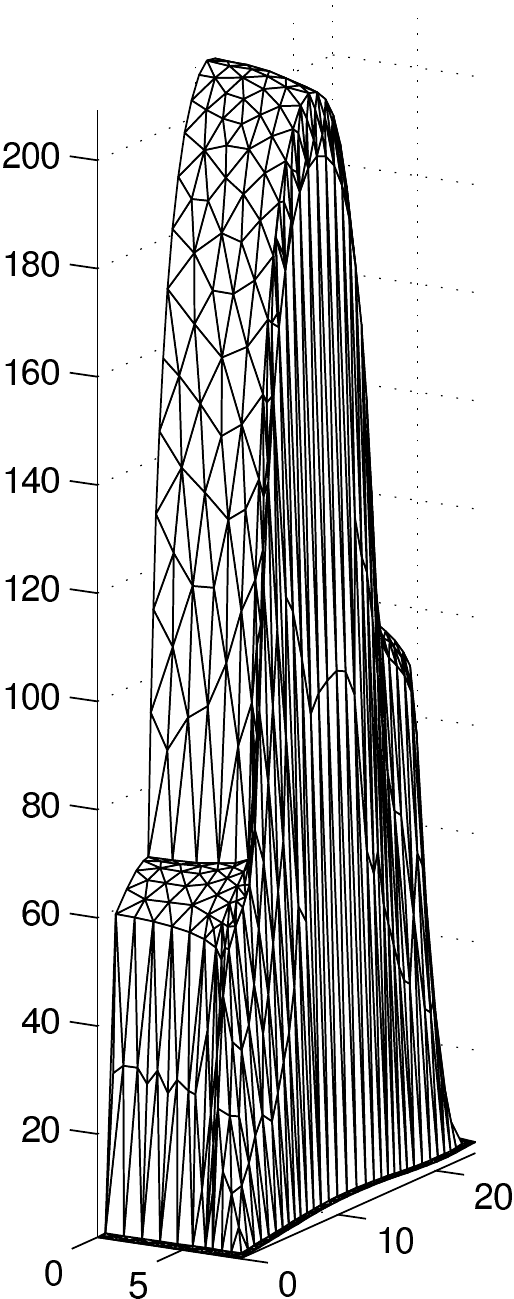}
   }%
   \caption{Heat conduction in a thermal battery: 
      (\protect\subref{fig:battery_geometry}) device geometry, 
      (\protect\subref{fig:battery_solution_}) contour plot, and 
      (\protect\subref{fig:battery_surface}) surface plot of a linear finite 
         element solution.}\label{fig:battery}
\end{figure}

The analytical solution for this problem is unavailable. The geometry and the contour and surface plots of a finite element approximation are given in \cref{fig:battery}.

We compare the quadratic least squares Hessian recovery and the full error estimator. 
For this example we found that three steps of the symmetrical Gauß--Seidel method were not sufficient for a full mesh adaptation and increased the number to seven, which proved to be enough to achieve at least a comparable error estimate as the one obtained with quadratic least squares Hessian recovery. 

\begin{figure}[t]
   \subcaptionbox{Interface edges are included in the mesh.\label{fig:battery_convergence_with}}{
      \includegraphics[width=0.47\textwidth,clip]{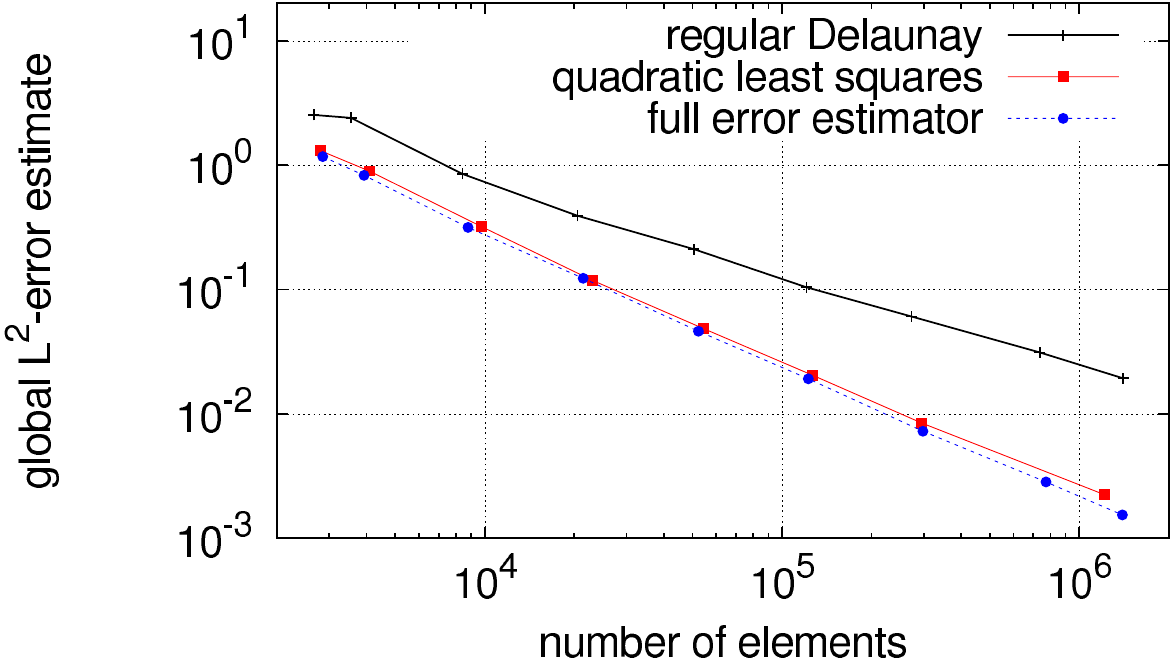}
   }%
   \hfill{}%
   \subcaptionbox{Interface edges are not included in the mesh.\label{fig:battery_convergence_without}}{
      \includegraphics[width=0.47\textwidth,clip]{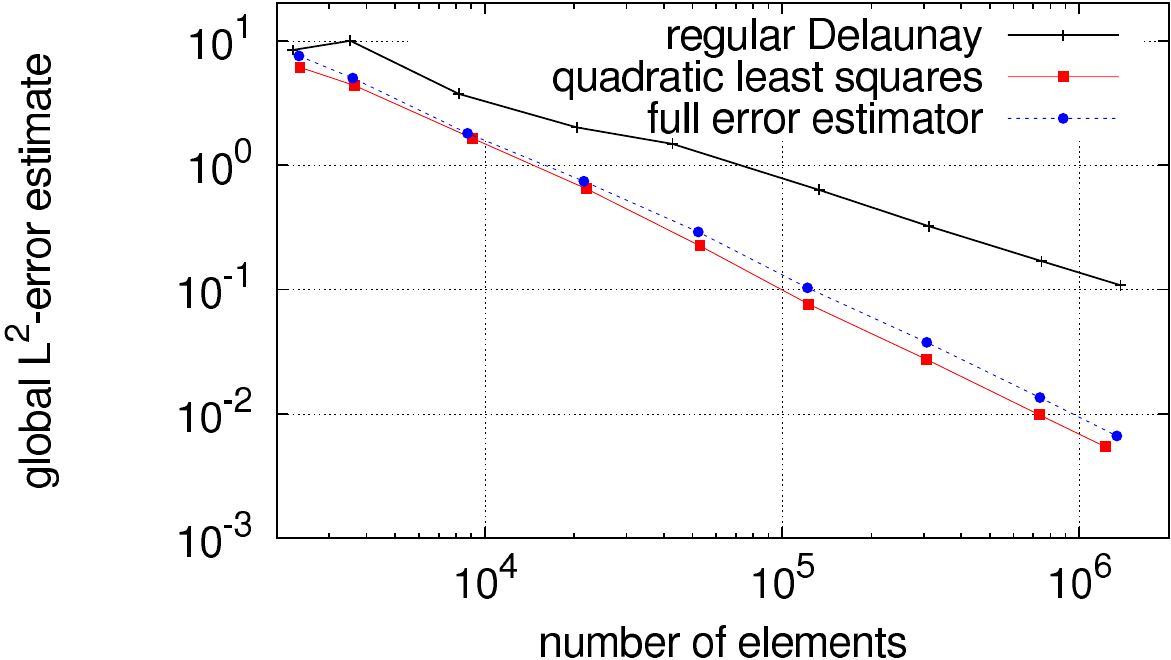}
   }%
   \caption{Heat conduction in a thermal battery: a comparison of the error for adaptive finite element solutions obtained on meshes (\protect\subref{fig:battery_convergence_with}) with and  (\protect\subref{fig:battery_convergence_without}) without the interfaces being present in the mesh.%
   }\label{fig:battery_convergence_plots}
\end{figure}

\cref{fig:battery_convergence_plots} shows global error estimates (obtained by solving exactly the approximate error problem $(E_h)$) for finite element solutions on adaptive meshes controlled by the full error estimate or Hessian recovery and having all or no predefined interface edges.
(The interface consists of edges when a mesh has all predefined interface edges.)

\begin{figure}[p] \centering
   \subcaptionbox{Quadratic least squares Hessian recovery:
      \num{3499} vertices and \num{6781} triangles,
      maximum aspect ratio $39.2$,
   error estimate $\|z_h\|_{L^2} = \num{4.7e-1}$.\label{fig:battery_qls}}[0.47\linewidth]{
      \includegraphics[width=0.41\textwidth,height=1.16\textwidth,clip]{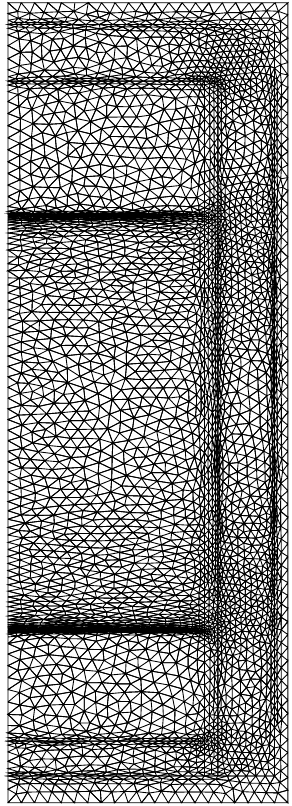}
   }\qquad
   \subcaptionbox{Full error estimator: 
      \num{3493} vertices and \num{6750} triangles,
      maximum aspect ratio $54.8$,
   error estimate $\|z_h\|_{L^2} = \num{4.2e-1}$.\label{fig:battery_hb_sgs07}}[0.47\linewidth]{
      \includegraphics[width=0.41\textwidth,height=1.16\textwidth,clip]{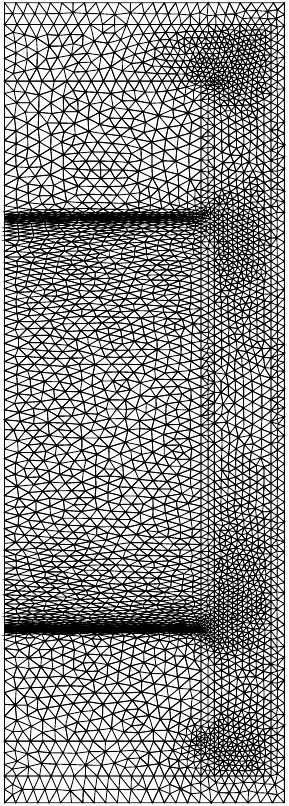}
   }
   \caption{Heat conduction in a thermal battery: adaptive meshes 
      obtained with (\protect\subref{fig:battery_qls}) quadratic least squares 
      Hessian recovery and (\protect\subref{fig:battery_hb_sgs07}) full  error
      estimator.}\label{fig:battery_qls_hb}
\end{figure}

Typical adaptive meshes with predefined interface edges for both methods are shown in \Cref{fig:battery_qls_hb}.

The results are in good agreement with those in \cref{sec:jG}.
When the interface edges are not present in the mesh, both methods provide similar results. 
On the other hand, when the mesh contains all the information of the interface, the quadratic least squares Hessian recovery produces a mesh with strong element concentration near all internal interfaces (\cref{fig:battery_qls}), whereas the full error estimator leads to a mesh (cf.~\cref{fig:battery_hb_sgs07}) that has higher element concentration in the corners of the regions, has a proper element orientation near the interfaces between the regions $2$ and $3$, and is almost uniform in regions where the solution is nearly linear (cf.~\cref{fig:battery_surface} for the surface plot of a computed solution).

Meshes without predefined interface edges are quite similar to those in the example with discontinuous gradients (\cref{sec:jG,fig:jg_without_interface}).
The interfaces are recognized by the both methods and the obtained adaptive meshes are dense near the interfaces.

Once again, the numerical results for this example show that a recovery method can lead to over-concentration of elements.
The new method, on the other hand,  produces only necessary concentration and is also able to catch the directional information of the solution required for proper element alignment.
This example also demonstrates that the new method can be successfully used for problems with jumping coefficients and strong anisotropic features.

\section{Conclusions and comments}\label{sec:Conclusion}

In the previous sections, we have presented a mesh adaptation method based on hierarchical basis error estimates and shown that anisotropic mesh adaptation can be successfully controlled by a posteriori error estimators. 
Numerical results have shown that the new method is fully comparable in accuracy with commonly used Hessian recovery-based methods and can be more efficient for some examples by producing only necessary element concentration.

A key idea in the new approach is the use of the full hierarchical error estimator for reliable directional information of the solution.  
To avoid the expensive exact solution of the global error problem, we employed only a few steps of the symmetric Gauß--Seidel iteration for the efficient solution of the resulting linear system.
Numerical results have shown that this is sufficient for obtaining an approximation to the error good enough for the purpose of mesh adaptation.

\section*{Acknowledgments} The work was partially supported by the German Research Foundation (DFG) under grants SFB568/3 and SPP1276 (MetStroem) and by the National Science Foundation (USA) under grants DMS-0410545 and DMS-0712935.

The authors are grateful to the anonymous referees for their valuable comments.

\newpage{}
\bibliographystyle{elsarticle-num}
\bibliography{HuKaLa10}
\end{document}